\documentclass{article}
\usepackage{amsmath}
\usepackage{amsfonts}
\usepackage{amssymb}
\usepackage{graphicx}
\usepackage{amsthm}
\usepackage{ascmac}
\usepackage{comment}
\usepackage{ulem}
\usepackage{color}
\usepackage{bm}

\theoremstyle{plain}
\newtheorem{theo}{Theorem}[section]

\newtheorem{defn}{Definition}[section]
\newtheorem{rem}{Remark}[section]
\newtheorem{lem}{Lemma}[section]
\newtheorem{prop}{Proposition}[section]
\newtheorem{ass}{Assumption}[section]

\theoremstyle{definition}
\newtheorem{ex}{Example}[section]

\newcommand{\ncd}{\newcommand}

\newcommand{\beq}{\begin{eqnarray}}
\newcommand{\eeq}{\end{eqnarray}}

\newcommand{\beqs}{\begin{eqnarray*}}
	\newcommand{\eeqs}{\end{eqnarray*}}

\ncd{\benu}{\begin{enumerate}}
	\ncd{\eenu}{\end{enumerate}}

\ncd{\btheo}{\begin{theo}}
	\ncd{\etheo}{\end{theo}}

\ncd{\bprop}{\begin{prop}}
	\ncd{\eprop}{\end{prop}}

\ncd{\blem}{\begin{lem}}
	\ncd{\elem}{\end{lem}}

\ncd{\bdef}{\begin{defn}}
	\ncd{\eddef}{\end{defn}}

\ncd{\bpf}{\begin{proof}}
	\ncd{\epf}{\end{proof}}

\ncd{\ds}{\displaystyle}

\ncd{\z}{{\mathbb Z}}
\ncd{\n}{{\mathbb N}}
\ncd{\R}{{\mathbb R}}
\ncd{\rd}{{\mathbb R}^d}
\ncd{\Ex}{{\mathbb E}}

\ncd{\calB}{{\mathcal B}}
\ncd{\calF}{{\mathcal F}}
\ncd{\calG}{{\mathcal G}}
\ncd{\calH}{{\mathcal H}}
\ncd{\calK}{{\mathcal K}}

\ncd{\bfB}{{\mathbf B}}
\ncd{\bfE}{{\mathbf E}}
\ncd{\bfX}{{\mathbf X}}
\ncd{\bfP}{{\mathbf P}}

\ncd{\prob}{{\mathbb P}}
\ncd{\intd}{\mathrm{d}}
\ncd{\wt}{\widetilde}
\ncd{\ol}{\overline}

\ncd{\indi}{1\hspace{-.2em}{\rm l}}
\ncd{\supp}{{\rm supp}}

\allowdisplaybreaks[4]
\makeatletter

\@addtoreset{equation}{section}
\makeatother
\begin{document}
\title{Limiting distributions for particles near the frontier of 
spatially inhomogeneous branching Brownian motions}

\author{Yasuhito Nishimori\thanks{Department of General Education,
National Institute of Technology, Anan College,
Anan, Tokushima, 774-0017, Japan; 
\texttt{nishimori@anan-nct.ac.jp}}}
\maketitle

\begin{abstract}
Our purpose in this paper is to determine the limiting distribution and the evolution rate of 
particles near the frontier of branching Brownian motions. 
Here the branching rate is given by a Kato class measure with compact support in Euclidean space. 
Our investigation focuses on the two dimensional case. 
\end{abstract}

\section{Introduction}
We consider branching Brownian motions with splitting on a compact set in $\rd$. 
The maximal displacement $L_t$ is the maximum of the Euclidean norm of particles at time $t$. 
Then, $L_t$ grows linearly and the growth rate is determined by $\lambda < 0$, 
the principal eigenvalue of the Schr\"odinger-type operator induced by the branching Brownian motion
(\cite{BH14}, \cite{E84} and \cite{S18}). 
Therefore, the frontier of particles lies around the boundary of a ball with a linear-growing radius centered at the origin. 
The aim of this research is to investigate asymptotic behaviors of the population size and 
distribution of particles near the frontier. 

We first explain the model of branching Brownian motions on $\rd$. 
A particle starts at $x \in \rd$ and moves according to the law of a standard Brownian motion $\{ B_t , t \ge 0\}$. 
Let $\mu $ be a branching rate measure, that is, a random branching time $T$ is distributed by 
\beq
	\bfP_x \left( T>t \mid B_s , s \ge 0 \right)
	= 
	e ^{-A_t ^{\mu}} .
	\label{eq:Ltas}
\eeq 
Here $\bfP_{x}$ is the law of the branching Brownian motion initiated from $x$ and 
$A_t ^{\mu}$ is the positive continuous additive functional (PCAF for short) 
which is in the Revuz correspondence with $\mu$. 
If $\mu$ is the Lebesgue measure, then $T$ has the exponential distribution with parameter one and 
the process is a spatially-homogeneous branching Brownian motion (HBBM for short); 
otherwise, the process is said to be inhomogeneous or catalytic. 
In particular, the branching rate depends on the trajectories of particles and no branch occurs outside the support of $\mu$.  

The branching Brownian motion on $\R$ has the rightmost particle at each time $t$ and its distance 
from the origin corresponds to $L_t$. McKean (\cite{Mc75} and \cite{Mc76}) proved that for a binary-HBBM, 
the distribution function of $L_t -R(t)$ converges to a unique solution of the F-KPP traveling-wave equation at speed 
$\sqrt{2}$ if $R(t) \sim \sqrt{2} t$. 
He also remarked on the logarithmic correction of $R(t)$. 
Bramson \cite{B78} revealed it and Mallein \cite{M15} extended it to $d \ge 2$ and gave an estimate of 
the tail probability $\bfP _x (L_t > R(t))$. 
On the other hand, Erickson \cite{E84}, Bocharov and Harris \cite{BH14}
for $d=1$, and Shiozawa \cite{S18} for $d \ge 1$ proved the linear growth of $L_t$ ((\ref{Ltlim}) below) 
for various inhomogeneous-BBMs. 
In addition, Lalley and Sellke \cite{LS88}, Bocharov and Harris \cite{BH16}, and 
Nishimori and Shiozawa \cite{ShioandN} proved that when 
\[
	R(t) = \sqrt{\frac{-\lambda}{2}} t + \frac{d-1}{2\sqrt{-2\lambda}} \log t , 
\] 
the tail distribution of $L_t - R(t)$ converges to the Gumbel distribution 
of which the parameter is mixed by the limit of a martingale ((\ref{mart}) below). 
Bocharov and Harris \cite[Theorem 2]{BH14}, and Shiozawa \cite[Theorem 2.8]{S18} proved the following 
for the number of particles outside a ball with radius $\delta t$ centered at the origin:  
if $\delta > \sqrt{-\lambda/2}$ ({\it supercritical}), 
then it converges to $0$ as $t \to \infty$. 
On the other hand, if $0 < \delta < \sqrt{-\lambda/2}$ ({\it subcritical}), 
then it increases exponentially, almost surely. 
Recently, Bocharov \cite{Bocha2020} showed that  
a distribution of particles near the frontier converges to the Poisson one 
for a BBM with a single-point catalyst at the origin, that is, 
in the case where $d=1$, $\mu = \beta \delta_0$. 
For the Dirac measure at the origin $\delta _0$, the PCAF $A_t ^{\beta \delta _0}$ is $\beta \ell _t$, 
where $\ell _t$ is  the local time at the origin of the Brownian motion. 
In this case, we see that $\sqrt{-\lambda/2}=\beta /2$ and thus
the frontier is either $- \beta t /2$ or $\beta t /2$. 
In this paper, we extend the results of Bocharov \cite{Bocha2020} to a BBM with $d \ge 1$ so that 
the branching rate measure is a Kato class measure $\mu$ with compact support in $\rd$. 
To do so, we develop the moment calculus of the population and the uniformly asymptotic behavior of the 
Feynman-Kac semigroup induced by $e^{A_t^\mu}$ as in \cite{ShioandN}. 
We need the scaling factor to determine the limiting distribution and the asymptotic growth 
of particles near the frontier. 
The significance of our research is to reveal that the scaling factor $R(t)$ below has the log correction 
$a(t)$ for $d = 2$ and, particularly, the critical $\delta$. 

For convenience, we use an annulus to explain our results. 
Let $R (t) = \delta t + a(t)$, where $\delta > 0$ and $a(t) \ge 0$ with $a(t) = o(t)$ as $t \to \infty$. 
Let us denote by $B_0 (R)$ a ball of radius $R$, centered at the origin. 
We write $A(R)=B_0 (R+r_2) \setminus B_0(R+r_1)$ for a $d$-dimensional annulus. 
When $\delta = \sqrt{-\lambda/2}$ and $t$ is large, $A(R(t)) $ is a domain near the frontier. 
For $A \subset \rd$, let $N_t ^A$ denote the number of particles in $A$ at time $t$. 
We claim the following (see Theorems \ref{theorem1}--\ref{theorem3}):
\begin{enumerate}
\item[{\rm (i)}] 
(Supercritical case) For all $d \ge 1$, if $\delta \in (\sqrt{-\lambda /2}, \sqrt{-2\lambda} )$,  
then $\bfP _x ( N_t ^{A(R(t))} > 0 )$ converges exponentially to zero. 

\item[{\rm (ii)}] (Subcritical case) For $d=1,2$ and $\delta \in (0,\sqrt{-\lambda/2})$, 
$N_t ^{A(R(t))}$ is increasing exponentially in probability. 

\item[{\rm (iii)}] (Critical case) 
For $d=1,2$ and $\delta = \sqrt{-\lambda/2}$, we take $a(t)$ as in (\ref{eq:a(t)}) below. 
Then, $\bfP _x ( N_t ^{A(R(t))} > 0 )$ converges polynomially to zero.  
Additionally, under some condition for $a(t)$, the distribution of $N_t ^{A(R(t))}$ converges to the Poisson one. 
Inside near the boundary of $B_0(R(t))$, the number of particles slowly grows. 
\end{enumerate}
As mentioned above, Bocharov studied the one-dimensional BBM with a single-point catalyst. 
In \cite[Proposition 2]{Bocha2020}, he showed (i), (ii) for $\delta \in (\beta /2 , \beta)$
and determined the limiting distribution in (iii) for $\delta = \beta /2$.  
However, (i) and (ii) were not clear in the critical case. 
By attaching the log correction $a(t)$ in (\ref{eq:a(t)}) to the leading term $\sqrt{-\lambda /2} t$, 
we reveal the critical case as follows: 
Theorem \ref{theorem1} establishes a polynomial decay of $\bfP _x ( N_t ^{A(R(t))} > 0 )$ and 
Theorem \ref{theorem2} gives the non-exponential growth order of the number of frontier particles. 
These results correspond to (i) and (ii), respectively. 
As a result, we can conclude that these asymptotic behaviors have the phase transition between 
$\delta= \sqrt{-\lambda/2}$ and $\delta \in (\sqrt{-\lambda /2}, \sqrt{-2\lambda})$. 
These are new results not contained in \cite{Bocha2020}. 

Our proofs are similar to those introduced by Bocharov \cite{Bocha2020}. 
He computed the asymptotic behavior of the distribution of particles near the frontier 
by using the first and second moments of the number of particles. 
By the Many-to-One Lemma, the first order moment is represented by the Feynman-Kac functional. 
He computed it directly by using the joint distribution of Brownian motion and local time. 
We show it for more general cases and use an analytic method established in \cite{ShioandN}.  
By the Many-to-Two Lemma and the crucial estimate of the Feynman-Kac semigroup, 
we show that the second order moment is asymptotically the same as the first order moment. 

In Section 2, we introduce the notions of branching Brownian motions and present our results. 
In Section 3, we compute the first and second moments of the population size near the frontier 
by using long time asymptotic properties of Feynman-Kac semigroups. 
Section 4 is devoted to the proofs of our results.

Throughout this paper, the letters $c$ and $C$ (with subscript and superscript) 
denote finite positive constants which may vary from place to place. 
For positive functions $f(t)$ and $g(t)$ on $(0,\infty)$, we write 
$f(t) \lesssim g(t)$, $t \rightarrow \infty$ if positive constants $T$ and $c$ exist such that 
$f(t)\le c g(t)$ for all $t \geq T$.  
We also write $f(t)\sim  g(t)$, $t \rightarrow \infty$ if $f(t)/g(t)\rightarrow 1$ as $t \rightarrow \infty$. 
We will omit ``$t \to \infty$'' for short when the meaning is clear.  
\section{Frameworks and results} 
We use the same notation as in \cite{ShioandN}. 
\subsection{Notations and some facts}
Let $(\{ B_t \}_{t \ge 0}, \{P_x \}_{x \in \rd}, \{\calF _t \}_{t \ge 0})$ be the Brownian motion on $\rd$ 
and $p_t (x,y)$ its transition function, where $\{\calF _t\}$ is the minimal augmented admissible filtration. 
For $\alpha > 0$, the $\alpha$-resolvent density $G_\alpha (x,y)$ of the Brownian motion is given by 
\[
	G_\alpha (x,y) = \int _0 ^\infty e^{-\alpha t}p_t (x,y) \intd t, 
	\quad x,y \in \rd , \ t > 0 .
\]
\begin{defn} 
\begin{itemize}
\item[{\rm (i)}] A positive Radon measure $\nu$ on ${\mathbb R}^d$  
is in the Kato class {\rm (}$\nu\in {\cal K}$ in notation{\rm )} if 
\[
	\lim_{\alpha\rightarrow\infty}
	\sup_{x\in {\mathbb R}^d}\int_{{\mathbb R}^d}G_{\alpha}(x,y) \nu(\intd y)=0.
\] 
\item[{\rm (ii)}] For $\beta>0$, a measure $\nu \in {\cal K}$ is $\beta$-Green tight 
{\rm (}$\mu\in {\cal K}_{\infty}(\beta)$ in notation{\rm )} if
\[
	\lim_{R\rightarrow\infty} \sup_{x\in {\mathbb R}^d} 
	\int_{{|y|\geq R}}G_{\beta}(x,y) \nu(\intd y)=0.
\] 
When $d \geq 3$, $\nu \in {\cal K}$ belongs to ${\cal K}_{\infty}(0)$ if 
the equality above is valid for $\beta=0$.
\end{itemize}
\end{defn}
We know by \cite{T08} that ${\cal K}_{\infty}(\beta)$ is independent of $\beta > 0$. 
Any Kato class measure with compact support is $1$-Green tight by definition. 

For $\nu\in \calK$, let $A_t ^{\nu}$ 
be the positive continuous additive functional associated with $\nu$ 
under the Revuz correspondence (see \cite[p.401]{FOT11}). 
For a signed measure $\nu=\nu^{+}-\nu^{-} \in {\cal K} - \calK$, 
we define $A_t^{\nu}=A_t^{\nu^+}-A_t^{\nu^-}$. 
The Feynman-Kac semigroup $\{p_t^{\nu}\}_{t>0}$ is defined by 
\[
	p_t^{\nu}f(x):=E_x\left[e^{A_t^{\nu}}f(B_t)\right], 
	\quad f \in {\cal B}_b({\mathbb R}^d)\cap L^2({\mathbb R}^d),
\] 
where ${\cal B}_b({\mathbb R}^d)$ is the collection of 
all bounded Borel measurable functions on ${\mathbb R}^d$. 
By \cite[Theorem 6.1 (ii)]{ABM91},  $\{p_t^{\nu}\}_{t>0}$ is a strongly continuous symmetric 
semigroup on $L^2({\mathbb R}^d)$. 
The corresponding $L^2$-generator is  called  a Schr\"odinger-type operator ${\cal H}^{\nu}=-\Delta/2-\nu$. 
Since $\{p_t^{\nu}\}_{t>0}$ is extended to $L^p({\mathbb R}^d)$ 
for any $p\in [1,\infty]$ by \cite[Theorem 6.1 (i)]{ABM91}, 
we  use the same notation $\{p_t^{\nu}\}_{t>0}$ as the extended one. 
By \cite[Theorems 7.1]{ABM91}, 
$p_t ^\nu$ possesses a jointly continuous integral kernel $p_t^{\nu}(x,y)$ 
on $(0,\infty)\times {\mathbb R}^d\times {\mathbb R}^d$ such that 
$$ 
	p_t^{\nu}f(x)=\int_{{\mathbb R}^d}p_t^{\nu}(x,y)f(y)\,{\rm d}y, 
	\quad  f \in {\cal B}_b({\mathbb R}^d) .
$$

For $\nu = \nu^+ - \nu^{-}\in {\cal K}_{\infty}(1) - {\cal K}_{\infty}(1)$, 
let us denote by $\lambda (\nu)$ the bottom of the spectrum for ${\cal H}^{\nu}$:
\begin{equation*}
\lambda(\nu)=
\left. \inf\left\{\frac{1}{2}\int_{{\mathbb R}^d}|\nabla u|^2\,{\rm d}x
-\int_{{\mathbb R}^d}u^2\,{\rm d}\nu \ \right| \  
u\in C_0^{\infty}({\mathbb R}^d), \int_{{\mathbb R}^d}u^2\,{\rm d}x=1\right\},
\end{equation*}
where $C_0^{\infty}({\mathbb R}^d)$ is the collection of all smooth functions on ${\mathbb R}^d$ 
with compact support. 
If $\lambda(\nu)<0$, then $\lambda(\nu)$ is the principal eigenvalue of ${\cal H}^{\nu}$ 
(\cite[Lemma 4.3]{T03} or \cite[Theorem 2.8]{T08}) and $h$ is the corresponding eigenfunction. 
Then $h$ has a strictly positive, bounded and continuous version (\cite[Section 4]{T08}). 
We also write $h$ for this version with $L^2$-normalization $\|h\|_{L^2({\mathbb R}^d)}=1$. 
Hence for any $x\in {\mathbb R}^d$ and $t>0$,
\begin{equation*}
p_t^{\nu}h(x)=e^{-\lambda t}h(x).
\end{equation*}
We assume that both $\nu^{+}$ and $\nu^{-}$ are compactly supported in ${\mathbb R}^d$. 
By the proof of \cite[Theorem 5.2]{T08} or \cite[Appendix A.1]{S19}, 
there exist positive constants $c_1$ and $c_2$ such that  
\begin{equation}\label{eq:est-eigen}
\frac{c_1e^{-\sqrt{-2\lambda(\nu)}|x|}}{|x|^{(d-1)/2}}\leq h(x)
\leq \frac{c_2e^{-\sqrt{-2\lambda(\nu)}|x|}}{|x|^{(d-1)/2}} ,
\quad |x|\geq 1 .
\end{equation}

Let $\lambda_2(\nu)$ be the second bottom of the spectrum for ${\cal H}^{\nu}$:
\begin{equation*}
\lambda_2(\nu)
:=\left. 
\inf\left\{\frac{1}{2}\int_{{\mathbb R}^d}|\nabla u|^2\,{\rm d}x-\int_{{\mathbb R}^d}u^2\,{\rm d}\nu 
\ \right| \ u\in C_0^{\infty}({\mathbb R}^d), 
\int_{{\mathbb R}^d}u^2\,{\rm d}x=1, \int_{{\mathbb R}^d}uh\,{\rm d}x=0\right\}.
\end{equation*}
If $\lambda(\nu)<0$, then $\lambda(\nu)<\lambda_2(\nu)\leq 0$ 
because the essential spectrum of ${\cal H}^{\nu}$ is the interval $[0,\infty)$ 
by \cite[Theorem 3.1]{BEKS94} or \cite[Lemma 3.1]{Be04}. 

\subsection{Branching Brownian motions}\label{subsect:branching}
In this subsection, we introduce the branching Brownian motion
(see \cite{INW68-1,INW68-2,INW69} and \cite{S18,S19} for details). 
Let $\mu \in \calK$ be a branching rate and $\{ p_n (x) , n \ge 1 \}$ a branching mechanism, 
where 
\[
	0 \le p_n (x) \le 1 , \ n \ge 1 \ \text{ and } 
	\sum _{n =1} ^{\infty} p_n (x) = 1, \ x \in \rd .
\]
A random time $T$ has an exponential distribution
\[
	P _x \left( T > t \mid \calF _{\infty} \right) = e^{-A_t ^{\mu}}, \quad t > 0. 
\]
A Brownian particle starts at $x \in \rd$. After an exponential random time $T$, 
the particle splits into $n$ particles with probability $p_n (B_T)$. 
New ones are independent Brownian particles starting at $B_T$ and 
each one independently splits into some particles, the same as the first. 
The $n$ particles are represented by a point in the following configuration space $\bfX$. 
Let $(\rd)^{(0)} = \{ \Delta \}$ and $(\rd)^{(1)} = \rd$. For $n \ge 2$, we define the equivalent 
relation $\sim$ on $(\rd)^n = \underbrace{\rd \times \cdots \times \rd}_{n}$ as follows: 
for $\bm{x} ^n = (x^1,\dots , x^n)$ and $\bm{y} ^n = (y^1,\dots , y^n) \in (\rd ) ^n$, 
we write $\bm{x} \sim \bm{y}$ if there exists a permutation $\sigma $ on $\{1,2,\dots , n\}$ such that 
$y^i = x^{\sigma (i)}$ for any $i \in \{1,2\dots ,n \}$. If we define $(\rd)^{(n)}=(\rd)^n \slash \sim $ for 
$n \ge 2$ and ${\bf X} = \bigcup _{n=0}^{\infty} (\rd)^{(n)}$, then $n$ points in $\rd$ determine a point 
in $(\rd)^{(n)}$. The branching Brownian motion $(\{ \bfB _t \}_{t \ge 0}, 
\{ \bfP_{\bm{x}} \}_{\bm{x} \in \bfX}, \{ \calG _t \}_{t \ge 0})$ is an $\bfX$-valued Markov process. 
Abusing notation, we regard $x \in \rd$ in the same way as $\bm{x} \in (\rd)^{(1)}$ and 
write $\bfP _x$ for $x \in \rd$. That is, $(\{ \bfB _t \}_{t \ge 0} , \bfP _x , \{ \calG _t \}_{t \ge 0})$ is 
the branching Brownian motion such that a single particle starts from $x \in \rd$. 

Let $Z_t$ be the set of all particles and $\bfB_t ^u$ the position of $u \in Z_t$ at time $t$. 
For $f \in \calB _b (\rd )$, 
\[
	Z_t (f) := 
		\ds\sum _{u \in Z_t} f \left( \bfB _t ^u \right), \quad t \ge 0 .
\]
For $A \subset \rd$, we set $N_t ^A = Z_t (\indi _A)$, in particular, $N_t = Z_t (1)$. 
The random variable $N_t^A$ is the number of particles which stay on $A$ and 
$N_t$ is the total number of particles at time $t$.  
Similarly, we use $Z_t ^A$ to denote a set of particles on $A$ at time $t$.

Assume that $\nu$ is a Kato class measure with compact support in $\rd$ and $\lambda := \lambda (\nu) < 0$. 
Let $h$ be the eigenfunction of $\calH ^{\nu}$ corresponding to $\lambda$ and 
\beq
	M_t := e^{\lambda t} Z_t (h), \quad t \ge 0 .
\label{mart} 
\eeq 
By the same argument as in \cite[Lemma 3.4]{S08}, we see that $M_t$ is 
a square integrable non-negative $\bfP _x$-martingale. Therefore, the limit 
$M_{\infty} := \lim _{t \to \infty} M_t \in [0,\infty)$ exists $\bfP _x$-a.s. and 
$\bfP _x (M_{\infty} > 0) > 0$. In particular, $\bfP _x (M_{\infty} > 0) =1$ for $d=1,2$ 
by \cite[Remark 2.11]{S18}. 

Let
\[
	L_t = \max _{ u \in Z_t} \left| \bfB _t ^u \right| .
\]
By \cite[Corollary 2.9]{S18}, 
\beq
	\lim _{t \to \infty} \ds\frac{L_t}{t} = \sqrt{\ds\frac{-\lambda}{2}}, \quad 
	\bfP_x \left( \cdot \mid M_{\infty} > 0 \right) \text{-a.s.}
	\label{Ltlim}
\eeq 
For $d=1,2$, (\ref{Ltlim}) holds $\bfP _x$-a.s. 

Let us recall the Many-to-One and Many-to-Two Lemmas. 
Let
\[
	Q(x) = \sum _{n=1} ^{\infty} n p_n (x), \quad R(x) = \sum _{n = 2} ^{\infty} n (n-1) p_n (x) .
\]
For a measure $\mu$, $Q\mu$ and $R \mu$ denote the measure $Q(x) \mu (\intd x)$ and $R(x) \mu (\intd x)$, 
respectively. 
\begin{lem}[{\cite[Lemma 3.3]{S08} and \cite{ShioandN}}]\label{lem:many-to}
Let $\mu \in \calK$. 
\begin{enumerate}
\item[{\rm (i)}] If the measure $Q \mu$ also belongs to Kato class, then for any $f \in \calB _b (\rd)$, 
\[
	\bfE _x \left[ Z_t (f) \right] = E _x \left[ e^{A_t^{(Q-1)\mu}} f \left( B_t \right) \right].
\]

\item[{\rm (ii)}] If the measure $R \mu$ also belongs to Kato class, then for any $f,g \in \calB _b (\rd)$, 
\[
	\bfE _x \left[ Z_t (f) Z_t (g) \right] 
	= 
	E _x \left[ e^{A_t ^{(Q-1)\mu}} f \left( B_t \right) g \left( B_t \right) \right]
	+
	E _x \left[ \int _0 ^t e^{A_s^{(Q-1)\mu}} \bfE _{B_s} \left[ Z_t (f) \right] \bfE _{B_s} \left[ Z_t (g) \right]
	\intd A_s ^{R\mu} \right] .
\] 
\end{enumerate}
\end{lem}

\subsection{Results}\label{subsection:Results}
We will make the following assumptions: 
\begin{ass}\label{ass} 
\begin{enumerate}
\item[{\rm (i)}] $\mu$ is a Kato class measure with compact support in $\rd$. 

\item[{\rm (ii)}] $R\mu \in \calK$.  

\item[{\rm (iii)}] $\lambda :=\lambda ((Q-1)\mu)<0$. 
\end{enumerate}
\end{ass}
 
Let $r_1 ,r_2\in \R$ with $r_1 < r_2$ and $\Theta \subset S^{d-1}$, where $S^{d-1}$ ($S^0 =\{ -1,1 \}$) 
is a unit sphere. 
We fix $r_1,r_2$ and $\Theta$. 
For $R(t)= \delta t + a(t)$, $\delta > 0$, $a(t) \ge 0$ and $a(t) = o(t)$ as $t \to \infty$, we set
\[
	C (R(t)) = \left\{ s x ; s \in [R(t)+r_1 , R(t)+r_2], x \in \Theta \right\} ,
\]
which represents a set near the frontier if $\delta = \sqrt{-\lambda /2}$. 
When $\Theta = S^{d-1}$, $C(R(t))$ is the annulus $A(R(t))$. 
Let $(R,\Theta)=\{ x \in \rd \mid |x|>R,x/|x| \in \Theta \}$, $R>0$, $\Theta \subset S^{d-1}$. 
From \cite[Remark 3.2]{ShioandN}, we see that 
\[
	\int _{(R,\Theta)} h(y) \intd y 
	\sim 
	c_{d,\lambda,\Theta} e^{-\sqrt{-2\lambda}R}
	R^{(d-1)/2}
	, \quad R \to \infty ,
\] 
where 
\beq
	c_{d,\lambda,\Theta} = \ds\frac{\sqrt{-2\lambda}^{(d-5)/2}}{(2 \pi)^{(d-1)/2}}
	\int _{\rd} \left( \int _{\Theta} e^{\sqrt{-2\lambda} \langle \theta ,z\rangle } \intd \theta \right) 
	h(z) \mu (\intd z ), 
	\quad c_{d}:=c_{d,\lambda, S^{d-1}} .
	\label{eq:c_d_s}
\eeq 
Let $\theta$ be the surface measure on $S^{d-1}$. If $d=1$, then $\theta = \delta _{-1} + \delta _1$, 
where $\delta _a$ is the Dirac measure on $a$.  
Thus we have 
\beq
	\int _{C(R(t))} h(y) \intd y
	\sim 
	c_* e ^{-\sqrt{-2\lambda} R(t)}
	R(t)^{(d-1)/2} , 
	\label{eq:c_*}	 
\eeq 
where 
$c_* =c_{d,\lambda ,\Theta} \left( e^{- \sqrt{-2\lambda} r_1} - e^{- \sqrt{-2\lambda} r_2} \right)$. 

We define 
\beq
	a_d (t) = \ds\frac{d-1}{2 \sqrt{-2\lambda}} \log (t \vee 1) + \gamma (t) .
	\label{eq:a(t)}
\eeq 
\begin{theo}\label{theorem1} 
Let Assumption \ref{ass} hold and let $\delta \in [\sqrt{-\lambda/2}, \sqrt{-2\lambda})$. 
When $\delta = \sqrt{-\lambda/2}$, we additionally take $a(t)=a_d(t)$ and $\gamma (t) \to \infty$. 
Then, 
	$$ 
	\lim _{t \to \infty} e^{\lambda t + \sqrt{-2\lambda} R(t)} R(t)^{-(d-1)/2}
	\bfP _x \left( N_t ^{C (R(t))} > 0 \right)
	= c_* h(x), \quad \text{ for all } x \in \rd . 
	$$ 
\end{theo} 
By the Paley-Zygmund inequality ((\ref{eq:Paley-Zygmund}) below) and 
the moment calculations of $N_t ^{C(R(t))}$, 
we show that 
\[
	\bfP _x (N_t ^{C(R(t))} > 0) \sim \bfE_x \left[ N_t ^{C(R(t))} \right] 
	\sim c_* h(x) e^{-\lambda t - \sqrt{-2\lambda}R(t)} R(t)^{(d-1)/2} .
\]   
In the case of $\delta = \sqrt{-\lambda /2}$, we particularly need (\ref{eq:a(t)}) and $\gamma (t) \to \infty$ to 
give the asymptotic lower bound of $\bfP _x (N_t ^{C(R(t))} > 0)$ 
(see Proposition \ref{lem:2ndorder-estimate} and (\ref{eq:bottleneck}) below), and then
	\beq 
	\lim _{t \to \infty} \left( \dfrac{-\lambda}{2} \right) ^{-(d-1)/4}
	e^{\sqrt{-2\lambda} \gamma (t)}
	\bfP _x \left( N_t ^{C (R(t))} > 0 \right)
	= c_* h(x), \quad \text{ for all } x \in \rd .
	\label{eq:criticalgamma}
	\eeq 

We next show the convergence in probability of the normalization of $N_t ^{C(R(t))}$. 
\begin{theo}\label{theorem2}
	Let Assumption \ref{ass} hold. Suppose $d=1$, $2$. 
	\begin{enumerate}
	\item[{\rm (i)}] If $\delta \in (0 , \sqrt{-\lambda/2})$, then for any $x \in \rd$, 
	\[
		e^{\lambda t + \sqrt{-2 \lambda} R(t)} R(t)^{-(d-1)/2} N_{t} ^{C (R(t))} \to c_* M_{\infty}, 
		\quad t \to \infty \quad \text{in probability } \bfP _x .
	\]
	\item[{\rm (ii)}] 
	For $\delta = \sqrt{-\lambda/2}$, 
	we set $a(t) = a_d (t)$ and assume that $\gamma (t) \to - \infty$ and $\gamma (t) = o (\log t)$. 
	Then 
	\[
		(-\lambda /2) ^{-(d-1)/4} e^{\sqrt{-2\lambda}\gamma (t)}N_t ^{C(R(t))} \to c_* M_\infty , 
		\quad t \to \infty \quad \text{in probability } \bfP_{x}. 
	\]
	\end{enumerate}
\end{theo}
Let $r_{1,i} , r_{2,i} \in \R$ with $r_{1,i} < r_{2,i}$ and $\Theta _i \subset S^{d-1}$. 
For each $i =1,2,\dots, n$, we write
\[
	C^i (R(t)) = \left\{ s x ; s \in [R(t) + r_{1,i} , R(t)+r_{2,i} ] , x \in \Theta _i \right\} .
\]
We always suppose that $\{ C^i (R(t)) ; i=1,2,\dots n \}$ are disjoint sets. 
Let $c_* (i)$ be the constant in (\ref{eq:c_*}) for $r_1=r_{1,i}, r_2=r_{2,i}$, $\Theta=\Theta _i$, and 
\[
	\wt{c} := \sum _{i=1}^n c_* (i) .
\]
We finally show that, if $d=1,2$, then the number of particles near the frontier converges in distribution 
to the Poisson-like distribution. 

\begin{theo}\label{theorem3}
	Let Assumption \ref{ass} hold and $\delta = \sqrt{-\lambda /2}$. 
	For $d=1,2$, we set $a(t)=a_d (t)$ and $\gamma (t) \ge 0$. 
	For any $ x \in \rd$ and $k _i \in \n \cup \{0 \}$, $i=1,2,\dots ,n$, 
	if $\gamma (t) \to \gamma \in [0,\infty)$, then 
\[
	\lim _{t \to \infty} 
	\bfP _x \left( \bigcap _{i=1} ^n \left\{ N_t ^{C^i(R(t))} = k_i  \right\} \right)
	= 
	\bfE _x \left[ 
		\exp \left( - \wt{c} M_\infty e^{ - \sqrt{-2 \lambda} \gamma } \right) 
		\ds\prod _{i=1} ^n 
		\dfrac{
		\left( c_* (i) M_{\infty} e^{- \sqrt{-2\lambda} \gamma } \right)^{k_i}
		}{k _i !}
	\right] 
\]
Here $k=k_1+\dots +k_n$. If $\gamma (t) \to \infty$, then it is the degenerate 
distribution on the origin. 
\end{theo}
\begin{rem}
If $\delta = \sqrt{-\lambda /2}$ and $n=1$, 
then the distribution is the Poisson one with parameter $c_*M_\infty e^{- \sqrt{-2\lambda} \gamma }$. 
In addition, if $k=0$, then it is the Gumbel one. 
This result was shown by \cite[Theorem 1.2]{BH16} and \cite[Theorem 2.4]{ShioandN}. 
\end{rem}
\begin{ex}\label{ex:1}
	We consider the one-dimensional binary-BBM. 
	Here $p_2 (x ) \equiv 1 $ and $Q(x) =R(x) \equiv 2$, then  
	the corresponding Schr\"odinger-type operator is $- \Delta / 2 - \mu$, 
	for the branching rate measure $\mu$. 
	
	(i) We assume that the branching rate measure is $\beta \delta_0$, $\beta >0$. 
	The binary-BBM splits only at the origin. 
	Since $\lambda = -\beta ^2 / 2$, Assumption \ref{ass} holds. 
	The corresponding 
	$L^2$-normalized eigenfunction is $h (x) = \sqrt{\beta} e^{-\beta |x|}$. 
	We see by (\ref{Ltlim}) the frontier particle sites near $\beta t /2$ or $- \beta t /2$ at time $t$. 
	Thus almost all particles are contained in $B_0(\beta t/2):=(-\beta t/2 , \beta t / 2)$ and
	the frontier particles hang around the boundary of $B_0(\beta t/2)$. 
	We define an $r$-neighborhood of the boundary of $B_0(R(t))$ as follows: 
	\[
		A (R(t),r) = (-R(t)-r, -R(t)] \cup [R(t),R(t)+r) ,
	\] 
	and consider the frontier particles on it for each $R(t)$. 
	Now the constant $c_*$ becomes $2 \sqrt{\beta}^{-1}  ( e^{\beta r}  - e^{-\beta r } )$. 
	When $d=1$, we write $R(t) = \delta t + \gamma (t)$. 
	\begin{enumerate}
	\item[(a)] If $\delta \in (0,\beta/2)$, then 
	$N_t ^{A(\delta t + o(t),r)}$ exponentially grows in probability, by Theorem \ref{theorem2}. 
	
	\item[(b)] Let $\delta \in ( \beta /2 , \beta)$. 
	For any $\gamma (t) = o(t)$, 
		\[
			e^{-\beta^2 t/2  + \beta R(t)} 
			\sim e^{(-\beta^2 /2 + \beta \delta ) t} \to \infty \quad \text{as } t \to \infty . 
		\]
	Thus, $\bfP _x (N_t ^{A(\delta t + \gamma (t),r)} > 0) $ exponentially converges to zero, 
	by Theorem \ref{theorem1}. 
	In other words, the particle rarely appears on $A(\delta t + \gamma (t),r)$, for large $t$.  
	
	\item[(c)]  
	Let $\delta = \beta /2$. 
	(1) We take $\gamma_0 (t)=(\gamma /\beta) \log t$, for any $\gamma >0$. 
	By (\ref{eq:criticalgamma}), 
	\[
		\lim _{t \to \infty} t^\gamma \bfP _x \left( N_t ^{A(\beta t/2 + \gamma_0(t),r)} > 0 \right) 
		= 
		c_* \sqrt{\beta} e^{-\sqrt{\beta} |x|} .
	\]
	Namely, we can find the particles on $A(\beta t/2 + \gamma _0(t),r)$ with high probability. 
	However, the number of frontier particles slowly grows. 
	Indeed, for $\gamma _1(t) = - \beta ^{-1} \sqrt{\log t}$,
	we see form Theorem \ref{theorem2} that 
	\[
		e^{- \sqrt{\log t}} N _t ^{A(\beta t / 2 +\gamma _1 (t) , r)}
		\to c_* M_{\infty} .
	\]

	(2) We next consider the limiting distribution of $N_t ^{A(\beta t/2 + \gamma (t),r)}$, 
	for $\gamma (t) \equiv \gamma$ or $\gamma (t) \to \gamma \in [0, \infty)$. 
	Noting that domain $A(R(t),r)$ is constructed from the two symmetrical sets with respect to the origin,  
	we see from Theorem \ref{theorem3} that $N_t ^{A(\beta t/2 + \gamma (t),r)}$ asymptotically follows 
	the Poisson one with parameter $c_* M_{\infty} e^{- \beta \gamma /2}$. 
	On the other hand, for any $\gamma (t) \to \infty$, 
	the distribution converges to a degenerate distribution on zero, for example $\gamma(t) = \gamma _0 (t)$. 
	\end{enumerate}
	Thus we can guess that the frontier particles stay on $A ( \beta t / 2, o(\log t) )$ and 
	the number of these particles follows the Poisson distribution. 
	In \cite{Bocha2020}, Bocharov proved (a), (b) and (c)-(2). 
	Our crucial result is (c)-(1). 
	
	Using the case of $k=0$, we have the limiting distribution for the maximal displacement $L_t$. 
	Let $R(t) = \beta t/2 + \gamma$ for any $\gamma >0$. 
	We set 
	$
		{\mathcal A} (R(t)) = (- \infty , R(t)] \cup [R(t) , \infty ) 
	$. 
	Then $c_* = 2 /\sqrt{\beta}$ and 
	\beqs
		\lim _{t \to \infty} \bfP_{x} \left( L_t \ge \dfrac{\beta}{2} t + \gamma \right) 
		&=& 
		\lim _{t \to \infty} \bfP _x \left( N_t ^{{\mathcal A}(R(t))} > 0 \right)
		=
		1 -\lim _{t \to \infty} \bfP _x \left( N_t ^{{\mathcal A}(R(t))} = 0 \right)
		\\
		&=& 
		1 - \bfE _x \left[ \exp \left( - c_* M_\infty  e^{- \beta \gamma} \right) \right] .
	\eeqs 
	Namely, the limiting distribution of $L_t$ is the Gumbel one. 
	This result was proved by Bocharov and Harris in \cite{BH16}. 
	
	(ii) We assume that the branching rate measure is $\mu = \delta _a + \delta _0$. 
	Then the binary-BBM splits only on $\{ 0,a \}$. 
	This model does not fall under the category of \cite{Bocha2020}. 
	Since $p=q=1$ in \cite[Example 3.10 (ii)]{S19}, we see that $\lambda < 0$. 
	Thus Assumption \ref{ass} is fulfilled and Theorems \ref{theorem1}--\ref{theorem3} hold. 
	The choice of $\gamma(t)$ is similar to (i). 
\end{ex}

\begin{ex}
We consider the $d$-dimensional ($d \ge 2$) binary-BBM splitting on surface 
$S_R^{d-1} = \{ x \in \R ^d \mid |x| =R \}$.  
The BBM has the Schr\"odinger-type operator $-\Delta /2 - \beta \theta _R$, 
where $\theta _R$ is the surface measure on $S_R^{d-1}$. 
We see the condition from \cite[Example 2.14]{S18} and references therein that 
$\lambda := \lambda (\beta \theta _R)< 0$ for some $\beta , R>0$. 
Thus Theorems \ref{theorem1}--\ref{theorem3} hold for the appropriate $d$, $\delta$ and $a(t)$. 
Similarly to Example \ref{ex:1} (i), we consider the frontier particles for $d=2$. 
Here, $A(R(t),r)$ is the annulus centered at the origin with width $r$: 
\[
	A(R(t), r) = B _0 (R(t)+ r) \setminus B _0 (R(t)).
\] 
When $d=2$ and $\delta = \sqrt{-\lambda /2}$, we take  
\[
	R_c (t) = \sqrt{\dfrac{-\lambda}{2}} t + \dfrac{1}{2\sqrt{-2\lambda}} \log (t \vee 1) + \gamma (t) .
\]
\begin{enumerate}
\item[$\bullet$]  If $\delta \in (\sqrt{-\lambda /2}, \sqrt{-2\lambda})$, 
then $\bfP _x ( N_t ^{A(\delta t + a(t), r)} > 0 )$ exponentially converges to zero for any $a(t) = o(t)$, 
by Theorem \ref{theorem1}.

\item[$\bullet$] 
Taking $\delta = \sqrt{-\lambda / 2}$ and 
$\gamma (t) = \gamma \sqrt{-2\lambda} ^{-1} \log (t \vee 1) $, 
$\gamma>0$, we see from (\ref{eq:criticalgamma}) that
		\[
			\lim _{t \to \infty} \left( \dfrac{-\lambda}{2} \right)^{-1/4}
			t^{\gamma} \bfP _x \left( N_t ^{A(R_c(t) ,r)} > 0 \right) 
			= c_* h(x) .
		\]
By Theorem \ref{theorem2}, 
for $\gamma (t) = - \sqrt{-2\lambda}^{-1} \log \log t$, 
\[
	\left( \dfrac{-\lambda}{2} \right)^{-1/4} (\log t )^{-1} N_t ^{A (R_c (t), r)} \to c_* M_\infty, 
	\quad \text{in probability}.
\]
Thus, the frontier particles stay around $S_{R_c(t)}^1$ with high probability, but, 
the number of frontier particles slowly increases. 
By Theorem \ref{theorem3}, if $\gamma (t) \to \gamma \in [0,\infty)$ or $\gamma (t) \equiv \gamma$, 
then $N_t ^{A(R_c (t),r)}$ asymptotically follows the Poisson distribution. 
On the other hand, the distribution of $N_t ^{A(R(t),r)}$ degenerates, 
if $\gamma (t) \to \infty$, even if $\gamma (t) = o(\log t)$. 
\end{enumerate}
Hence, we infer that the frontier particles are distributed on 
\[
	A \left( \sqrt{\dfrac{-\lambda}{2}}t + \dfrac{1}{2\sqrt{-2\lambda}} \log t, o(\log t) \right) 
\]
and the number on this set has the Poisson distribution. 
\end{ex}

\section{Estimates of the first and second moments}\label{sec3}
Let $R(t) = \delta t + a(t)$, where $\delta \in (0 , \sqrt{-2\lambda})$ and $a(t) = o(t)$ as $t \to \infty$. 
In this section, we calculate the moments of $N_t ^{C(R(t))}$ by using Lemma \ref{lem:many-to} and 
the Feynman-Kac semigroup $p_t ^{\nu}$. We introduce a kernel  $q_t ^{\nu}$ to estimate $p_t ^{\nu}$:  
\beq
	q_t ^{\nu} (x,y):= p_t ^{\nu} (x,y) - p_t (x,y)  - e^{-\lambda t} h(x) h(y) , \ x ,y \in \rd , \ t \ge 0.
	\label{eq:q_t}
\eeq
In what follows, we always assume the following: 
\begin{enumerate}
\item[$\bullet$] 
	for some $\alpha \in (0,1)$, $ 0 \le s (t) <\alpha t$ for any $t \ge 0$ 
	and $s(t) = o(t)$ as $t \to \infty$; 
	
\item[$\bullet$] $b(t)\ge 0$ for any $t \ge 0$ and $b(t)=o(t)$ as $t \to \infty$; 

\item[$\bullet$] $x:[0,\infty) \to \rd$ satisfies $|x(t)| \le b(t)$ for all large $t > 0$. 

\end{enumerate}
\subsection{The first moment}
The aim of this subsection is to calculate the first moment of $N_t^{C(R(t))}$. 
The main claim is provided by Proposition \ref{lem:orderof1moment}. 
We concisely explain the outline of the proof  of Proposition \ref{lem:orderof1moment} and 
show details in Lemma \ref{lem:estP_0-1}--\ref{lem:K_3} and Proposition \ref{prop:PropSN3-3}. 
\begin{prop}\label{lem:orderof1moment}
	Let $\mu^+$ and $\mu^-$ be Kato class measures with compact support in $\rd$ and 
	$\lambda := \lambda ((Q-1)\mu)< 0$. 
	Then there exist $T >0$ and $\theta_{\pm} (t)$ such that 
	\[
		\theta _- (t)
		\le 
		\ds\frac{\bfE _{x(t)} \left[ N_{t-s(t)} ^{C(R(t))}\right]}
		{c_* h(x(t)) e^{-\lambda (t-s(t)) -\sqrt{-2\lambda} R(t)} 
		R(t)^{(d-1)/2}}
		\le 
		\theta _+ (t)
	\] 
	for all $t \ge T$. Here $\theta _{\pm} (t) \to 1$ as $t \to \infty$.
\end{prop}
\bpf
This follows by the same method as in \cite[(3.14)]{ShioandN}. 
By Lemma \ref{lem:many-to} and (\ref{eq:q_t}), 
\beq
	&&
	\bfE _{x(t)} \left[ N_{t-s(t)} ^{C(R(t))}\right] 
	= 
	E_{x(t)} \left[ e^{A_{t-s(t)} ^{(Q-1)\mu}} \ ; \ B_{t-s(t)} \in C(R(t)) \right]
	\nonumber 
	\\ 
	&=& 
	P_{x(t)} \left( B_{t-s(t)} \in C(R(t)) \right) + e^{-\lambda (t-s(t))} h(x(t)) \int _{C(R(t))} h(y) \intd y 
	+ \int _{C(R(t))} q_{t-s(t)} ^{(Q-1)\mu} (x(t),y) \intd y 
	\nonumber
	\\ 
	&=:&
	{\rm (I)}  + {\rm (II)} + {\rm (III)}.
	\label{eq:20200106-1}
\eeq
In Lemma \ref{lem:estP_0-1} (i) and Proposition \ref{prop:PropSN3-3} (ii) below,  
we will show that there exist $C >0$ and $T>0$ such that for all $t \ge T$, 
\beq
	{\rm (I)} , \ |{\rm (III)}| \le h(x(t)) e^{-C t} e^{-\lambda (t-s(t)) - \sqrt{-2\lambda}R(t)} .
	\label{eq:20200107-3}
\eeq 
Then we have by (\ref{eq:c_*}), 
$$
	\frac{{\rm (I)} + |{\rm (III)}|}{{\rm (II)}}
	\sim 
	\frac{{\rm (I)} + |{\rm (III)}|}{c_* h(x(t)) e^{-\lambda(t-s(t))-\sqrt{-2\lambda}R(t)}
	R(t)^{(d-1)/2}} ,
$$ 
and we see from (\ref{eq:20200107-3}) that the right-hand side is bounded above by $C' e^{-C t}$. 
Therefore, 
	\beqs 
		&& 
		\ds\frac{\bfE _{x(t)} \left[ N_{t-s(t)} ^{C(R(t))}\right]}
		{c_* h(x(t)) e^{-\lambda (t-s(t)) - \sqrt{-2\lambda}R(t)}
		R(t)^{(d-1)/2}}
		\\
		&\le &
		\ds\frac{\text{(II)}}
		{c_* h(x(t)) e^{-\lambda (t-s(t)) - \sqrt{-2\lambda}R(t)}
		R(t)^{(d-1)/2}}
		\left( 1 + \frac{{\rm (I)} + {\rm |(III)|}}{{\rm (II)}} \right) 
		\\
		&=& 
		\ds\frac{\int _{C(R(t))} h(y) \intd y}
		{c_* e^{- \sqrt{-2\lambda}R(t)} 
		R(t)^{(d-1)/2}} 
				\left( 1 + \frac{{\rm (I)} + {\rm |(III)|}}{{\rm (II)}} \right) 
		=: \theta _+ (t) , 
	\eeqs 
where $\theta _+ (t) \to 1$ as $t \to \infty$ by (\ref{eq:c_*}) and (\ref{eq:20200107-3}). 

Similarly to the above,   
	\[
		\left| \frac{{\rm (I)} - |{\rm (III)}|}{{\rm (II)}} \right|
		\lesssim C' e^{-Ct}
		\to 0 , \quad t \to \infty  , 
	\]
and 
	$$ 
		\ds\frac{\bfE _{x(t)} \left[ N_{t-s(t)} ^{C(R(t))}\right]}
		{c_* h(x(t)) e^{-\lambda (t-s(t)) - \sqrt{-2\lambda}R(t)}
		R(t)^{(d-1)/2}}
		\ge 
		\ds\frac{\int_{C(R(t))} h(y) \intd y}
				{c_* e^{ - \sqrt{-2\lambda}R(t)}
				R(t)^{(d-1)/2}} 
				\left( 1 + \frac{{\rm (I)} - | {\rm (III)}|}{{\rm (II)}} \right) 		
		=:\theta _- (t) ,
	$$ 
where  $\theta _-  (t)\to 1$ as $t \to \infty$. 
\epf 

We give an estimate of (\ref{eq:20200107-3}). 
By Lemma \ref{lem:estP_0-1} (i) and 
\[
	P_{x(t)} \left( B_{t-s(t)} \in C( R(t) ) \right) 
	\le 
	P_{x(t)} \left( \left| B_{t-s(t)} \right| \ge R(t) \right) ,
\]
we have (\ref{eq:20200107-3}) for (I). 
Let $\lambda = \lambda (\nu)$ for a signed measure $\nu = \nu^+-\nu^-$. 
\begin{lem}\label{lem:estP_0-1}
Let $\nu ^+$ and $\nu^-$ be Kato class measures such that $\lambda < 0$ and the support of $\nu$ 
is contained in $\ol{B_0 (M)}$ for some $M>0$. 
Then the following three assertions hold. 
\begin{enumerate}
\item[{\rm (i)}] There exist $C>0$ and $T \ge 0$ such that for all $t \ge T$, 
\beq
	P_{x(t)} \left( \left| B_{t-s(t)} \right| \ge R(t) \right) 
	\le 
	h(x(t)) e^{-Ct} e^{-\lambda (t-s(t)) - \sqrt{-2\lambda} R(t)} . 
\label{eq:estP_0-core}
\eeq 
\item[{\rm (ii)}]  
There exist $C>0$ and $T \ge 1$ such that for all $x \in \ol{B_0(M)}$, $t \ge T$ and $w \in [0,t-1]$, 
\[
	P_{x} \left( \left| B_{t-w} \right| \ge R(t) \right) 
	\le 
	Ce^{-\lambda (t-w) - \sqrt{-2\lambda} R(t)} t^{(d-2)/2}. 
\]
\item[{\rm (iii)}] 
There exist $C>0$ and $T \ge 0$ such that for all $x \in \ol{B_0(M)}$, $t \ge T$ and $w \in [0,1]$, 
\[
	P_{x} \left( \left| B_{w} \right| \ge R(t) \right) 
	\le 
	C e^{-(R(t) - M)^2/2} t^{d-2} . 
\]
\end{enumerate}
\end{lem}
\bpf 
(i) We have by the spatial homogeneity and scaling property of Brownian motion, for any $x \in \ol{B_0(M)}$, 
\beqs
	P_x \left( \left| B_{t-s} \right| \ge R(t) \right) 
	\le 
	P_0 \left( \left| B_{t-s} \right| \ge R(t) - M \right) 
	= 
	P_0 \left( \left| B_1 \right| \ge \frac{R(t) - M }{\sqrt{t-s}} \right) .
%	\label{eq:Pspatial}
\eeqs 
Hence if $|x(t)| \le b(t) $, then 
\[
	P_{x(t)} \left( \left| B_{t-s(t)} \right| \ge R(t) \right) 
	\le 
	P_0 \left( \left| B_1 \right| \ge \frac{R(t) - b(t)}{\sqrt{t-s(t)}} \right) .
\] 
We set $R_b(t) := R(t) - b(t)$.
Since  
\beq
	\int _R ^{\infty} e^{-u^2/2} u^{d-1} \intd u \sim e^{-R^2/2} R^{d-2} , \quad R \to \infty 
	\label{eq:Gauss-int-asym}
\eeq 
and $R_b (t)/\sqrt{t-s(t)} \to \infty$ as $t \to \infty$, 
there exists $c>0$ such that for any large $t$, 
\beqs 
	P_0 \left( \left| B_1 \right| \ge \frac{R_b(t)}{\sqrt{t-s(t)}} \right) 
	&=& 
	\left( 2 \pi \right)^{-d/2} \omega _d
	\int _{ R_b(t)/\sqrt{t-s(t)}} ^{\infty} e^{-u^2/2} u^{d-1} \intd u 	
	\\
	&\le &  
	c \exp \left\{ - \frac{1}{2} \left( \frac{R_b (t)}{\sqrt{t-s(t)} }\right)^2 \right\}
	\left( \frac{R_b(t)}{\sqrt{t-s(t)} }\right)^{d-2} ,
\eeqs 
where $\omega_d$ is the area of the unit ball in $\rd$. 
By (\ref{eq:est-eigen}) and $|x(t)| \le b(t)$ for large $t$,  
\[
	1 
	= \frac{h(x(t))}{h(x(t))} 
	\le 
	c h(x(t)) e^{\sqrt{-2\lambda }|x(t)|} |x(t)| ^{(d-1)/2} 
	\le 
	c h(x(t)) e^{\sqrt{-2\lambda }b(t)} 
	b(t) ^{(d-1)/2} .
\]
Thus, 
\beqs  
	P_0 \left( \left| B_1 \right| \ge \frac{R_b(t)}{\sqrt{t-s(t)}} \right)  
	&\le &
	h(x(t)) 
	\exp \left\{ 
		\log c - \frac{1}{2} \left( \frac{R_b(t)}{\sqrt{t-s(t)} }\right)^2  
	\right. \\ && \left. 
		+ (d-2) \log \frac{R_b(t)}{\sqrt{t-s(t)}}
		+ \sqrt{-2\lambda} b(t)
		+ \frac{d-1}{2} \log b(t)
	\right\} 
	\\ 
	&=& 
	h(x(t)) \exp \left\{ - \frac{1}{2} \left( \frac{R_b(t)}{\sqrt{t-s(t)} }\right)^2 + o(t) \right\} . 
\eeqs
Then
$$ 
	- \frac{1}{2} \left( \frac{R_b(t)}{\sqrt{t-s(t)}} \right)^2
	=
	\lambda (t-s(t)) \left\{ 1 -  \frac{R_b(t)}{\sqrt{-2\lambda} (t-s(t))} \right\}^2 + \sqrt{-2\lambda} b(t) 
	- \sqrt{-2\lambda}R(t) -\lambda (t-s(t)) 	
$$
and 
\[
	\lim _{t \to \infty} \frac{1}{t} 
	\left[
		\lambda (t-s(t)) \left\{ 1 -  \frac{R_b(t)}{\sqrt{-2\lambda} (t-s(t))} \right\}^2 
		+ \sqrt{-2\lambda} b(t) + o(t)
	\right]
	= 
	\lambda \left( 1 - \frac{\delta}{\sqrt{-2\lambda}} \right) ^2
	<0 .
\]
Hence we can choose $C>0$ and large $T>0$ such that for any $t \ge T$,  
$$ 
	- \frac{1}{2} \left( \frac{R_b(t)}{\sqrt{t-s(t)}} \right)^2 + o(t)
	\le 
	-C t 
	- \sqrt{-2\lambda}R(t) -\lambda (t-s(t)) 	.
$$
This proves (i). 

(ii) For any $x \in \ol{B_0 (M)}$ and $w \in [0,t-1]$, 
$$
	P_{x} \left( \left| B_{t-w} \right| \ge R(t) \right)
	\le 
	P_{0} \left( \left| B_{1} \right| \ge \frac{R(t) - M}{\sqrt{t-w}} \right) 
	= 
	( 2 \pi )^{-d/2} \omega _d 
	\int^\infty _{R_M(t) / \sqrt{t-w}} e^{-r^2/2} r^{d-1} \intd r , 
$$
where $R_M(t):=R(t) - M$. 
Then for any $w \in [0,t-1]$,  
\[
	\frac{R_M(t)}{\sqrt{t-w}} \ge \frac{R_M(t)}{\sqrt{t}} \to \infty, \quad t \to \infty .
\]
Therefore, in the same way as (i), for any large $t$,
$$
	\int^\infty _{R_M(t)/\sqrt{t-w}} e^{-r^2/2} r^{d-1} \intd r
	\le 
	c \exp \left\{ - \frac{1}{2} \left( \frac{R_M(t)}{\sqrt{t-w}} \right)^2 \right\} 
	\left( \frac{R_M(t)}{\sqrt{t-w}} \right) ^{d-2} .
$$
Then 
\beqs
	- \frac{1}{2} \left( \frac{R_M(t)}{\sqrt{t-w}} \right)^2
	&=& 
	- \left[ 
		\left\{ \sqrt{-\lambda (t-w)} - \frac{R_M(t)}{\sqrt{2(t-w)}} \right\}^2 + \sqrt{-2\lambda} R_M(t)
	\right] 
	- \lambda (t-w)
	\\
	&\le &
	- \sqrt{-2\lambda} R_M(t) - \lambda (t-w) ,
\eeqs
and $( R_M(t) / \sqrt{t-w} ) ^{d-2} \le c t^{(d-2)/2}$. 
Thus, there is some large $T>0$ such that for any $t \ge T$ and $w \in [0,t-1]$, 
\[
	P_x \left( \left| B_{t-w} \right| \ge R(t) \right) 
	\le 
	C e ^{- \lambda (t-w) - \sqrt{-2\lambda} R(t)} t^{(d-2)/2} .
\]

(iii) For any $x \in \ol{B_0 (M)}$ and $w \in [0,1]$, 
\[
	P_x \left( \left| B_{w} \right| \ge R(t) \right) 
	\le 
	P_0 \left( \left| B_1 \right| \ge R(t) - M \right) 
	= 
	(2 \pi ) ^{-d/2} \omega _d \int _{R_M(t)} ^{\infty} e^{-r^2/2} r^{d-1} \intd r .
\]
There exists $T > 0$ such that for any $t \ge T$, 
\[
	\int _{R_M(t)} ^{\infty} e^{-r^2/2} r^{d-1} \intd r 
	\le 
	c_1 e^{-R_M ^2 (t)/2} R ^{d-2} _M(t)
	\le 
	c_2 e ^{-( R(t) - M )^2 / 2} t^{d-2}
\]
and the proof is complete. 
\epf 

We estimate (III). 
By the same argument as \cite[(3.19)]{ShioandN}, 
\beqs
	\int _{C(R)} q_t^{\nu}(x,y) \intd y
	&=& 
	\int _0 ^1 \left( \int _{\rd} p_s ^{\nu} (x,z) P_z (B_{t-s} \in C(R)) \nu (\intd z) \right) \intd s
	\\ 
	&&
	+
	\int _1 ^t \left\{ \int _{\rd}\left( p_s ^{\nu} (x,z) -e^{-\lambda s} h(x)h(z)\right) 
	P_z (B_{t-s} \in C(R) ) \nu (\intd z) \right\} \intd s 
	\\ 
	&&
	-
	h(x) e^{-\lambda t} \int _{t-1} ^{\infty} e^{\lambda s} 
		\left( \int _{\rd}h(z) P_z (B_s \in C(R) ) \nu (\intd z) \right) \intd s
	\\
	&=:& 
	K_1 (x,t,C(R))+ K_2 (x,t,C(R))- K_3 (x,t,C(R)).
\eeqs 
We write $K_i (x,t,R) = K_i(x,t,C(R))$. The task is now to estimate the following: 
\[
	|{\rm (III)}|
	\le 
	|K_1 (x(t),t-s(t), R(t))| + |K_2 (x(t),t-s(t), R(t))| + | K_3 (x(t),t-s(t),R(t)) | . 
\]
We divide the estimate into a sequence of lemmas for each $K_i$. 
\begin{lem}\label{lem:K_1}
Under the same setting as in Lemma \ref{lem:estP_0-1}, the following three assertions hold. 
	\begin{enumerate}
	\item[{\rm (i)}] {\rm (\cite[Lemma 3.5]{ShioandN})}
	There exists $C>0$ such that for all $x \in \rd$, $t \ge 1$ and $R>M$, 
	\[ 
		\left| \int _0 ^1 \left( \int _{\rd} p_s ^{\nu} (x,z) 
			P_z \left( \left| B_{t-s} \right| \ge R \right) \nu (\intd z) \right) \intd s \right|  
		\le C h(x) P_0 \left( \left| B_t \right| > R -M \right) .
	\] 
	\item[{\rm (ii)}] There exist $C>0$ and $T>1$ such that for all $x \in \rd$, $t \ge T$ and $w \in [0,t-1]$, 
	\[
		| K_1 (x, t-w,R(t)) | \le C h(x) e^{-\lambda (t-w) -\sqrt{-2\lambda} R(t)} t^{(d-2)/2} . 
	\]
	\item[{\rm (iii)}] There exist $C>0$ and $T>0$ such that for all $t \ge T$, 
	\[
		| K_1 (x(t), t-s(t),R(t)) | \le h(x(t))e^{-Ct} e^{-\lambda (t-s(t)) -\sqrt{-2\lambda} R(t)}.
	\] 
	\end{enumerate}
\end{lem}
For a signed measure $\nu = \nu ^+ - \nu ^-$, we write $|\nu| = \nu^+ + \nu^-$. 
\bpf 
	(ii) By Lemma \ref{lem:estP_0-1} (ii), 
	there exists $T>1$ such that for any $t \ge T$, $R(t) + r_1 > M$ and 
	\[
		P_0 (|B_{t-w}| \ge R(t)+r_1 -M) \le C e^{-\lambda (t-w) - \sqrt{-2\lambda}R(t)} 
		t^{(d-2)/2} , \quad w \in [0,t-1] .
	\]
	Thus, we see from (i) that
	\beqs
		\left| K_1 (x , t-w , R(t)) \right|
		&\le &
		\left| 
			\int _0 ^1 \int _{\rd} p_s ^{\nu} (x,z) 
			P_z \left( \left| B_{t-w-s} \right| \ge R (t) + r _1 \right) \nu (\intd z)  \intd s  
		\right| 
		\\
		&\le &
		c_1 h(x) P_0 (|B_{t-w}| \ge R(t)+r_1 -M)
		\\ 
		&\le &
		c_2 h(x)  e^{-\lambda (t-w) - \sqrt{-2\lambda}R(t)} t^{(d-2)/2} .
	\eeqs 
	
	(iii) follows by (\ref{eq:estP_0-core}) and the same argument of (ii). 
\epf 
For $c>0$, we set  
\[
	I_c (t, R)
	=
	\begin{cases}
	e^{c t-\sqrt{2c} R} R ^{(d-1)/2} , & \lambda _2 < 0 , \\
	e^{c t-\sqrt{2c} R} R ^{(d-1)/2} \wedge  t P_0 \left( \left| B_t \right| > R-M \right) ,
	& \lambda _2 = 0 .
	\end{cases} 	
\]
In \cite[Lemma 3.6]{ShioandN}, the following lemma was shown. 
\begin{lem}\label{lem:SN3-6}
Under the same setting as in Lemma \ref{lem:estP_0-1}, the following three assertions  hold. 
	\begin{enumerate}
	\item[{\rm (i)}] 
	For any $c>0$ with $c \ge -\lambda _2$, there exists $C >0$ such that for all $t \ge 1$ and $R>2 M$, 
	\[
		\int _1 ^t \left( \int _{\rd} e^{-\lambda _2 s} 
		P_z \left( \left| B_{t-s} \right| > R \right) |\nu| (\intd z) \right) \intd s
		\le 
		C e^{ct - \sqrt{2c}R} R^{(d-1)/2} .
	\]
	\item[{\rm (ii)}]
	If $\lambda _2 =0$, then there exists $C >0$ such that for all $t \ge 1$ and $R>M$, 
		\[
			\int _1 ^t \left( \int _{\rd} 
			P_z \left( \left| B_{t-s} \right| > R \right) |\nu| (\intd z) \right) \intd s
			\le 
			C t P_0 \left( \left| B_t \right| > R-M \right) .
		\]
	\item[{\rm (iii)}] For any $c >0$ with $c \ge -\lambda_2$, there exists $C >0$ 
	such that for all $x \in \rd$, $t \ge 1$ and $R>M$, 
	\[
		\int _1 ^t \left\{
			\int _{\rd} \left| p_s ^{\nu} (x,z) - e ^{-\lambda s} h(x) h(z) \right| 
			P_z \left( \left| B_{t-s} \right| > R \right) | \nu | (\intd z)
		\right\} \intd s 
		\le C I_c (t,R) .
	\]
	\end{enumerate}
\end{lem}
Lemma \ref{lem:SN3-6} yields upper estimates of $K_2$. 
\begin{lem}\label{lem:K2}
Under the same setting as in Lemma \ref{lem:estP_0-1}, the following assertions hold. 
\begin{enumerate}
\item [{\rm (i)}]
For any $c > 0$ with $c \ge -\lambda_2$, 
there exists $C>0$ such that for any $x \in \rd$, $t \ge 1$ and $R> M -r_1$,
\beq
	\left| K_2 (x,t,R) \right| \le C I_c (t,R+r_1) .
	\label{eq:K2-1}
\eeq 
In particular, there exists $T \ge 1$ such that for any $x \in \rd$, $t \ge T$ and $w \in [0,t-1]$, 
	$$
	\left| K_2 (x,t-w,R(t)) \right| \le C e^{-\lambda (t-w) - \sqrt{-2\lambda}R(t)} t^{(d-1)/2}.
	$$ 
\item[{\rm (ii)}] 
There exist $C>0$ and $T \ge 1$ such that for any $t \ge T$, 
\[
	\left| K_2 (x(t), t-s(t) , R(t)) \right| \le h(x(t)) e^{-Ct} e^{-\lambda (t-s(t)) - \sqrt{-2\lambda} R(t)}.
\]
\end{enumerate}
\end{lem}
\bpf
(i) (\ref{eq:K2-1}) follows by Lemma \ref{lem:SN3-6} (iii). 
Substituting $c=-\lambda$ in (\ref{eq:K2-1}), we have  
\beqs
	\left| K_2 (x,t-w,R(t)) \right| 
	&\le &
		C I_{-\lambda} (t-w,R(t)+r_1)
	\le 
	C e^{-\lambda (t-w) - \sqrt{-2\lambda}(R(t)+r_1)} (R(t)+r_1) ^{(d-1)/2} 
	\\ 
	&\le & 
	C'
	e^{-\lambda (t-w) - \sqrt{-2\lambda} R(t)} t ^{(d-1)/2} , 
\eeqs 
for any large $t$. 

(ii) Let $\alpha \in (0,1-\delta / \sqrt{-2\lambda})$, which is the constant at the beginning of Section \ref{sec3}. 
Fix a constant $c_0$ with 
	\[
			\left\{
				(-\lambda_2) \vee \left( \frac{\sqrt{2}\delta}{1-\alpha} - \sqrt{-\lambda} \right)^2 
			\right\}
			< c_0 < -\lambda .
		\]
	Since $s(t) < \alpha t$, $R(t) >M-r_1$ and $t -s(t) \ge 1$ for any large $t$, 
	(\ref{eq:K2-1}) becomes 
	\[
		\left| K_2 (x(t), t-s(t) , R(t)) \right| 
		\le 
		c_1 I_{c_0} (t-s(t), R(t)+r_1) .
	\] 
	By (\ref{eq:est-eigen}), 
	\[
		c_1 I_{c_0} (t-s(t), R(t)+r_1)
		\le 
		c_2 h(x(t)) e^{\sqrt{-2\lambda}|x(t)|} |x(t)|^{(d-1)/2} I_{c_0} (t-s(t), R(t)+r_1) .
	\] 
	Since
	\[ 
		c_2 e^{\sqrt{-2\lambda}|x(t)|} |x(t)|^{(d-1)/2}
		\le 
		c_2 e^{\sqrt{-2\lambda}b(t)} 
		b(t)^{(d-1)/2}
		=
		e^{o(t)} ,
	\]  
	it suffices to show that there exists $c>0$ such that for any large $t$, 
	\beq
		e^{o(t)}  I_{c_0} (t-s(t), R(t)+r_1) \le e^{-C t} e^{-\lambda (t-s(t))-\sqrt{-2\lambda}R(t)} .
		\label{eq:20200215-1}
	\eeq 
	Note that 
	\[
		e^{o(t)}
		I_{c_0} (t-s(t), R(t)+r_1) 
		\le
		e^{o(t)}
		e^{c_0 (t-s(t)) - \sqrt{2c _0} ( R(t) + r_1 ) } ( R(t) + r_1 ) ^{(d-1)/2}
		=
		e^{ c_0 (t-s(t)) - \sqrt{2c_0} R(t) + o(t) } .
	\]
	Then, 
	\beqs
		&& 
		c_0(t-s(t)) - \sqrt{2 c_0} R(t) + o(t)
		\\
		&=&
		\left\{ 
			( c_0 +\lambda ) \left(1-\frac{s(t)}{t} \right) 
			- \frac{( \sqrt{2 c_0} - \sqrt{-2\lambda} ) R (t)}{t} + o(1)
		\right\} t 
		- \lambda (t-s(t)) - \sqrt{-2\lambda} R (t) 
	\eeqs 
	and 
	\beqs 
		&& 
		\lim _{t \to \infty} 
		\left\{ 
					( c_0 +\lambda ) \left(1-\frac{s(t)}{t} \right) 
					- \frac{( \sqrt{2 c_0} - \sqrt{-2\lambda} ) R (t)}{t} + o(1)
		\right\} 
		\\ 
		&=&
		c_0 +\lambda 
		- 
		( \sqrt{c_0} - \sqrt{-\lambda} ) \sqrt{2} \delta 
		< 
		(\sqrt{c_0} - \sqrt{-\lambda}) \ds\frac{\sqrt{2}\delta \alpha}{1-\alpha} < 0 .
	\eeqs 
	Therefore, there exists $C>0$ such that for any large $t$, 
	\[
		( c_0 +\lambda ) \left(1-\frac{s(t)}{t} \right) 
		- \frac{( \sqrt{2 c_0} - \sqrt{-2\lambda} ) R (t)}{t} + o(1)
		< - C
	\]
	and (\ref{eq:20200215-1}) is proved. 
\epf
Let us set
\[
	J (t,R) = e^{-\lambda t - \sqrt{-2\lambda}R} R^{(d-1)/2} 
	\int _{(\sqrt{-2\lambda} t-R)/\sqrt{2t}} ^{\infty} e^{-v^2} \intd v .
\]
\begin{lem}[{\cite[Lemma 3.7]{ShioandN}}]\label{lem:SN3-7}
Under the same setting as in Lemma \ref{lem:estP_0-1}, 
there exists $C>0$ such that for any $x \in \rd$, $t \ge 1$ and $R>2M$, 
\[
	e^{-\lambda t} h(x) \int _{t-1} ^{\infty} e^{\lambda s} 
	\left( \int _{\rd} h(z) P_z \left( \left| B_s \right| > R \right) | \nu | (\intd z)\right) \intd s
	\le 
	C h(x) \left( P_0 \left( \left| B_t \right| >R-M \right) + J (t,R) \right) .
\]
\end{lem}

\begin{lem}\label{lem:K_3}
Under the same setting as in Lemma \ref{lem:estP_0-1}, the following assertions hold. 
\begin{enumerate}
\item[{\rm (i)}] 
There exists $C>0$ such that for any $x \in \rd$, $t \ge 1$ and $R>0$ with $R +r_1 > 2 M$, 
\[
	| K_3 (x , t,R) | \le C h(x) \left( P_0 (|B_{t}| > R+r_1-M) + J (t,R+r_1) \right) . 
\]
In particular, there exist $C'>0$ and $T\ge 1$ such that for any $x \in \rd$, $t \ge T$ and $w \in [0,t-1]$, 
\[
	| K_3 (x , t-w,R(t)) | \le C' h(x) 
	e^{-\lambda (t-w) - \sqrt{-2\lambda}R(t)} \left( t^{(d-2)/2} \vee t^{(d-1)/2} \right) .
\] 
\item[{\rm (ii)}] There exist $C>0$ and $T>0$ such that for any $t \ge T$,   
\[
	| K_3 (x(t) , t-s(t) ,R(t)) | \le h(x(t)) e^{-Ct} e^{-\lambda (t-s(t)) - \sqrt{-2\lambda}R(t)} .
\]
\end{enumerate}
\end{lem}
\bpf
	(i) For $t \ge 1$ and $R+r_1 > 2 M$. By Lemma \ref{lem:SN3-7}, for any $x \in \rd$, 
	\beqs
		\left|
			K_3 (x,t,R) 
		\right|
		&\le & 
			e^{-\lambda t} h(x) \int _{t-1} ^{\infty} e^{\lambda s} 
			\left( \int _{\rd}h(z) P_z (| B_s | > R+r_1) | \nu| (\intd z) \right) \intd s
		\\
		&\le & 
		Ch(x) \left( P_0 \left( \left| B_{t} \right| > R + r_1 - M \right) + J (t,R+r_1) \right) .
	\eeqs
	In particular, taking any large $t$ with $R(t) +r_1 > 2 M$, we see that for any $w\in [0,t-1]$, 
	\[
		\left|
			K_3 (x,t-w,R(t)) 
		\right|	
		\le 
		C h(x) \left( P_0 \left( \left| B_{t-w} \right| > R (t)+r_1 - M \right) + J (t-w,R(t)+r_1) \right) .
	\]
	Lemma \ref{lem:estP_0-1} (ii) yields  
	\[
		P_0 \left( \left| B_{t-w} \right| > R (t)+r_1 - M \right)
		\le 
		c e^{-\lambda(t-w) - \sqrt{-2\lambda}R(t)} t^{(d-2)/2}.
	\]
Then
\beqs 
	J(t-w, R(t)+r_1) 
	\le 
	e^{-\lambda (t-w) - \sqrt{-2\lambda}( R(t) + r_1)} (R(t) + r_1)^{(d-1)/2} 
	\le 
	c' e^{-\lambda (t-w) - \sqrt{-2\lambda} R(t)} t^{(d-1)/2}.
\eeqs 
	It follows that 
	\[
		\left| K_3 (x,t-w,R(t)) \right|	
		\le 
		c'' h(x) e^{-\lambda (t-w) - \sqrt{-2\lambda}R(t)} 
		\left( t^{(d-2)/2} \vee t^{(d-1)/2} \right) .
	\]  
	
	(ii) Note that $t-s(t) \ge 1$ and $R(t)+r_1 > 2M$ for any large $t$. 
	By (i), 
	\[ 
		\left| K_3 (x(t),t-s(t),R(t)) \right|
		\le  
		c_1 h(x(t)) \left( P_0 \left( \left| B_{t-s(t)} \right| > R (t) +r_1 - M \right) + J (t-s(t),R(t)+r_1) \right) .
	\] 
	Let us show that the right-hand side is bounded above by  
	\beq
		h (x(t)) e^{-C t} e^{-\lambda (t-s(t)) - \sqrt{-2\lambda}R(t)} .
		\label{eq:20200304-3} 
	\eeq 
	Lemma \ref{lem:estP_0-1} (i) implies that for large $t$, 
	\beq 
		c_1 P_0 \left( \left| B_{t-s(t)} \right| > R (t) +r_1 - M \right)
		&\le &
		c_1 \| h \|_{\infty}
		e^{-c_2t} e^{-\lambda (t-s(t)) -\sqrt{-2\lambda}( R(t) +r_1 -M)} 
		\nonumber 
		\\ 
		&\le &
		e^{-c_3 t} e^{-\lambda (t-s(t)) -\sqrt{-2\lambda} R(t) } ,
		\label{eq:20200214-3}
	\eeq
	where
	$c_1 \| h \| _{\infty} e^{-c_2 t} e^{-\sqrt{-2\lambda}(r_1 -M)} \le e^{-c_3 t}$.
	On the other hand, 	
	\beq
	\begin{aligned}
		c_1J (t-s(t),R(t)+r_1) 
		&= 
		c_1e^{-\lambda (t-s(t)) - \sqrt{-2\lambda}R(t)} (R(t)+r_1)^{(d-1)/2}
		\\ 
		& \quad \times 
		\int _{\{ \sqrt{-2\lambda} (t-s(t))-(R(t)+r_1)\}/\sqrt{2(t-s(t))}} ^{\infty} e^{-v^2} \intd v .
	\end{aligned}
	\label{eq:20210104-1}
	\eeq 
	By (\ref{eq:Gauss-int-asym}), 
	$
			\int _L^{\infty} e^{- v^2} \intd v 
			\sim 
			e^{- L^2} / (2 L)
	$ as $L \to \infty$. 
	Since for any $\delta \in (0 , \sqrt{-2\lambda})$, 
	\[
		L(t) :=\ds\frac{\sqrt{-2\lambda} (t-s(t))-(R(t)+r_1))}{\sqrt{2(t-s(t))}}
		= 
		\ds\frac{\sqrt{-2\lambda} - \delta + o(1)}{\sqrt{2 +o(1)}} \sqrt{t} \to \infty, \quad t \to \infty, 
	\]
	there exists some $c>0$ such that $L^2 (t) \ge ct$ for all sufficiently large $t$. 
	We thus have for some $C>0$, 
$$
		c_1 (R(t)+r_1)^{(d-1)/2}
		\int _{L(t)} ^{\infty} e^{-v^2} \intd v 
		\le 
		c' (R(t)+r_1)^{(d-1)/2}
		\frac{e^{-L^2(t)} }{L(t)}
		\le  
		e^{-C t} .
$$
	Therefore we have by (\ref{eq:20210104-1})
	\beq
		c_1 J(t-s(t), R(t)+r_1) \le 
		e^{-C t}e^{-\lambda (t-s(t)) - \sqrt{-2\lambda} R(t)} .
		\label{eq:20200630-1} 
	\eeq 	
	Combining (\ref{eq:20200214-3}) and (\ref{eq:20200630-1}), we obtain (\ref{eq:20200304-3}). 
\epf 
From Lemmas \ref{lem:K_1}, \ref{lem:K2} and \ref{lem:K_3}, we have the following proposition. 
In particular, when $\nu = (Q-1)\mu$, this gives an upper estimate of (III) in the proof of 
Proposition \ref{lem:orderof1moment}. 
\begin{prop}\label{prop:PropSN3-3}
Under the same setting as in Lemma \ref{lem:estP_0-1}, the following assertions hold.  
\begin{enumerate}
\item[{\rm (i)}] There exist $C>0$ and $T>0$ such that for all $x \in \rd$, $t \ge T$ and $w \in [0,t-1]$,  
\[
	\left| \int _{C(R(t))} q_{t-w} ^{\nu}(x,y) \intd y \right| 
	\le 
	C e^{-\lambda (t-w) - \sqrt{-2\lambda}R(t)} \left( t^{(d-2)/2} \vee t^{(d-1)/2} \right) .
\] 
\item[{\rm (ii)}]
There exist $C>0$ and $T>0$ such that for all $t \ge T$,  
\[
	\left| \int _{C(R(t))} q_{t-s(t)}^{\nu} (x(t),y) \intd y \right| 
	\le 
	h(x(t))e^{-C t} e^{-\lambda (t-s(t)) - \sqrt{-2\lambda}R(t)} .
	\label{eq:estq}
\]
\end{enumerate}
\end{prop}
\subsection{The second moment}\label{sec:2order}
We set 
\[
	C_d (t) = (d-1) \log (t \vee 1) , \quad  d \ge 1. 
\]
\begin{lem}\label{lem:2order}
Under the same setting as in Proposition \ref{lem:orderof1moment}, 
there exists $T \ge 1$ such that for any $t\ge T$, 
\beq
	\bfE _{x(t)} \left[ N_{t-s(t)} ^{C(R(t))} \right]
	\le 
	\bfE _{x(t)} \left[ \left( N_{t-s(t)} ^{C(R(t))} \right)^2 \right] 
	\le 
	\bfE _{x(t)} \left[ N_{t-s(t)} ^{C(R(t))} \right]
	+
	h(x(t))e^{- 2 \lambda (t-s(t))- 2 \sqrt{-2\lambda}R(t) + C_d (t)} .
	\label{eq:20210104-2}
\eeq 
\end{lem}
\bpf 
Since $N_{t-s(t)} ^{C(R(t))}$ is a non-negative integer, 
the first inequality of (\ref{eq:20210104-2}) holds. 

Let us denote by $\sigma _M$ the hitting time of some particles to $B_0 (M)$. 
Because $N_{t-s(t)} ^{C (R(t))} \in \{0,1\}$ on the event $\{ t - s(t) < \sigma _M \}$, 
\beq
	\bfE _{x(t)} \left[ \left( N_{t-s(t)} ^{C (R(t))} \right) ^2 \right] 
	=
	\bfE _{x(t)} \left[ \left( N_{t-s(t)} ^{C (R(t))} \right) ^2 ; t -s (t) \ge \sigma _M \right] 
	+ 
	\bfE _{x(t)} \left[ N_{t-s(t)} ^{C (R(t))} ; t -s (t) < \sigma _M \right] . 
	\label{eq:20200225-1} 
\eeq
It is sufficient to show that the first term of the right-hand side in (\ref{eq:20200225-1}) is bounded above by
\beq
	\bfE _{x(t)} \left[ N_{t-s(t)} ^{C(R(t))} ; t -s (t) \ge \sigma _M \right] 
	+ 
	C h(x(t)) e^{-2\lambda (t-s(t))-2 \sqrt{-2\lambda} R(t) + C_d (t)} , 
	\label{eq:20210104-3}
\eeq 
which gives the second inequality in  (\ref{eq:20210104-2}). 
Lemma \ref{lem:2ndorder} below shows that for large $t$, $u \le t-s(t)$, $s \in [0,t-s(t)-u]$ 
and $x \in \ol{B_0(M)}$, 
\beq 
	E _x \left[
		\int _0 ^{t-s(t)-u} e^{A_s ^{(Q-1)\mu}} 
		\bfE_{B_s} \left[ N_{t-s(t)-u-s} ^{C (R(t))} \right]^2 \intd A_s ^{R \mu}
	\right] 
	\le 
	C h (x) e^{-2\lambda (t-s(t)-u)-2\sqrt{-2\lambda} R(t) + C_d (t)} .
	\label{eq:20210105-1}
\eeq 
Before the initial particle hits $B_0 (M)$, the branching Brownian motion has no branch. 
Hence by the strong Markov property, Lemma \ref{lem:many-to} (ii) and (\ref{eq:20210105-1}), 
\beqs 
	&& 
	\bfE _{x(t)} \left[ \left( N_{t-s(t)} ^{C (R(t))} \right) ^2 ; t -s (t) \ge \sigma _M \right] 
	=
	\bfE _{x(t)} \left[ \left. \bfE _{\bfB _u} \left[ 
		\left( N_{t-s(t)-u} ^{C (R(t))} \right) ^2 \right] \right| _{u=\sigma _M} ; t -s (t) \ge \sigma _M
	\right] 
	\nonumber 
	\\
	&=&
	\bfE _{x(t)} \left[ N_{t-s(t)} ^{C (R(t))} ; t - s(t) \ge \sigma _M \right] 
	\nonumber
	\\ 
	&& \quad + 
	E_{x(t)} \left[ \left.  
	E _{B_u} \left[ 
				\int _0 ^{t-s(t)-u} e^{A_s ^{(Q-1)\mu}} 
				\bfE_{B_s} \left[ N_{t-s(t)-u-s} ^{C (R(t))} \right]^2 \intd A_s ^{R \mu}
			\right]
		\right| _{u=\sigma _M}	
	; t - s(t) \ge \sigma _M
	\right] 
	\nonumber
	\\
	&\le &
	\bfE _{x(t)} \left[ N_{t-s(t)} ^{C (R(t))} ; t - s(t) \ge \sigma _M \right]  
	\nonumber
	\\ 
	&& \quad + 
	C e^{-2\lambda (t-s(t))-2\sqrt{-2\lambda} R(t) + C_d (t)}
	E_{x(t)} \left[ 
		h (B_{\sigma _M}) e^{-2\lambda \sigma _M}
		; t - s(t) \ge \sigma _M
	\right] .
	\label{eq:20210104-4} 
\eeqs 
Since $M_t = e^{\lambda t} e^{A_t^\nu} h(B_t)$ is a martingale, the optional stopping theorem yields that 
\beqs
	& &
	E _{x(t)} \left[ h(B_{\sigma _M}) e^{2\lambda \sigma _M} ; t -s (t) \ge \sigma _M \right]
	= 
	E _{x(t)} \left[ h(B_{\sigma _M \wedge (t-s(t))}) e^{2\lambda \{ \sigma _M \wedge (t-s(t) )\} } 
	; t -s (t) \ge \sigma _M \right]
	\\
	&\le & 
	E _{x(t)} \left[ h(B_{\sigma _M \wedge (t-s(t))}) e^{\lambda \{ \sigma _M \wedge (t-s(t) )\} } \right]
	=
	h(x(t) ). 
\eeqs 
Therefore we have (\ref{eq:20210104-3}).
\epf 
\begin{lem}\label{lem:2ndorder} 
There exist $C>0$ and $T>0$ such that for any $x \in \ol{B_0(M)}$, $t \ge T$ and $u \in [0 , t-s(t)]$, 
	\[
		E _{x} \left[
			\int _0 ^{t-s(t)-u} e^{A_s ^{(Q-1)\mu}} 
			\bfE_{B_s} \left[ N_{t-s(t)-u-s} ^{C (R(t))} \right]^2 \intd A_s ^{R \mu}
		\right]
		\le 
		C h (x) e^{-2\lambda (t-s(t)-u)-2\sqrt{-2\lambda} R(t) + C_d (t)} .
	\]
\end{lem}
\bpf
According to Lemma \ref{lem:estP_0-1} (i) and (iii), 
we can take $t$ so large enough that
\[
	P_x \left( \left| B_{t-s(t)} \right| \ge R(t) \right) 
	\le 
	e^{-Ct} e^{-\lambda(t-s(t)) - \sqrt{-2\lambda}R(t)} 
\] 
and 
\beq
	P_0 ( | B_1 | \ge  R(t) + r_1 -M )
	\le 
	c e^{-(R(t)+r_1-M)^2/2} t^{d-2}.
	\label{T2}
\eeq
We first assume that $w_1 :=t-s(t) -u \le 1$. For any $s \in [0,w_1]$, 
\beq
	&& 
	\bfE_x \left[ N _{w_1-s} ^{C(R(t))} \right] 
	= 
	E_x \left[ e^{A_{w_1-s}^{(Q-1)\mu}} ; B_{w_1-s} \in C (R(t)) \right] 
	= 
	\int _{C (R(t))} p_{w_1-s} ^{(Q-1)\mu} (x,y) \intd y
	\nonumber
	\\
	&=& 
	\int _{C (R(t))} \left(p_{w_1-s} (x,y) + p_{w_1-s} ^{(Q-1)\mu} (x,y) - p_{w_1-s} (x,y) \right) \intd y
	\nonumber 
	\\
	&= & 
	P_x \left( B_{w_1-s} \in C (R(t)) \right) + 
	\int _{C (R(t))} \left(p_{w_1-s} ^{(Q-1)\mu} (x,y) - p_{w_1-s} (x,y) \right) \intd y
	\nonumber 
	\\
	&= & 
	P_x \left( B_{w_1-s} \in C (R(t)) \right) 
	+ 
	\int _{C (R(t))} \int _0 ^{w_1-s} \int _{\rd} 
	p_v ^{(Q-1)\mu} (x,z) 
	p_{w_1-s-v} 
	(z,y) \mu (\intd z) \intd v \intd y. 
	\label{eq:20200226-1}
\eeq
Here (\ref{eq:20200226-1}) follows from \cite[Lemma 3.1 (i)]{ShioandN}. 
The first and second terms of (\ref{eq:20200226-1}) are bounded, respectively, by   
$$
	P_0 \left( \left| B_1 \right| \ge R(t) + r_1 -M  \right) 
$$
and 
\beqs
	&& 
	\int _{C (R(t))} \int _0 ^{w_1-s} \int _{\rd} p_v ^{(Q-1)\mu} (x,z) 
	p_{w_1-s-v} (z,y) \mu (\intd z) \intd v \intd y 
	\\ 
	&\le & 
	\int _0 ^{w_1-s} \int _{\rd} p_v ^{(Q-1)\mu} (x,z) 
	P_z ( B_{w_1-s-v} \in C ( R(t) ) ) | \mu | (\intd z) \intd v 
	\\
	&\le &
	\int _0 ^{w_1-s} \int _{\rd} p_v ^{(Q-1)\mu} (x,z) 
	P_0 ( | B_{w_1-s-v} | \ge  R(t)+ r_1 -M ) | \mu | (\intd z) \intd v
	\\ 
	&\le &
	P_0 ( | B_1 | \ge  R(t) + r_1 -M )
	\int _0 ^1 \int _{\rd} p_v ^{(Q-1)\mu} (x,z) | \mu | (\intd z) \intd v 
	\\ 
	&\le & 
	c_3 P_0 ( | B_1 | \ge  R(t) + r_1 -M ) .
\eeqs
Thus, we see from (\ref{T2}) and (\ref{eq:20200226-1}) that  for any $x \in \ol{B_0 (M)}$, 
\beq
	\bfE_x \left[ N _{w_1-s} ^{C (R(t))} \right] 
	\le 
	c _1 P_0 ( | B_1 | \ge  R(t) + r_1 -M )
	\le 
	c_2 e^{-(R(t)+r_1-M)^2/2} t^{d-2} .
	\label{eq:20200226-4}
\eeq 
By \cite[Proposition 3.8]{CZ95}, 
\[
	\sup _{x \in \rd} E _{x} \left[ \left( A^{R|\mu|} _{t} \right)^2 \right]
	\le 
	\sup _{x \in \rd} E _{x} \left[ e^{A^{R|\mu|} _{t}} \right]
	\le 
	e^{c_1 + c_2 t}, \quad \text{ for all } t > 0, 
\]
and then for any large $t$, 
\beq
	&& 
	E _{x} \left[
		\int _0 ^{t-s(t)-u} e^{A_s ^{(Q-1)\mu}} \intd A_s ^{R \mu}
	\right]
	\le 
	E _{x} \left[
		\int _0 ^{t-s(t)-u} e^{A_s ^{(Q-1)\mu^+}} \intd A_s ^{R |\mu|}
	\right]
	\nonumber 
	\\ 
	&\le &
	E _{x} \left[
		e^{A_t ^{(Q-1)\mu^+}} \int _0 ^{t-s(t)-u} \intd A_s ^{R |\mu|}
	\right]	
	\le
	E _{x} \left[
		 e^{A_{t} ^{(Q-1)\mu^+}} A_{t} ^{R |\mu|}
	\right]	
	\nonumber 
	\\
	&\le & 
	E _{x} \left[ e^{2A_{t} ^{(Q-1)\mu^+}} \right]	^{1/2}
	E _{x} \left[ \left( A_{t} ^{R |\mu|} \right)^2 \right] ^{1/2}
	\le
	e^{c_3 t}.
	\label{eq:20201112-1}
\eeq
These imply that 
\beq
	&& 
	E _{x} \left[
		\int _0 ^{t-s(t)-u} e^{A_s ^{(Q-1)\mu}} 
		\bfE_{B_s} \left[ N_{t-s(t)-u-s} ^{C (R(t))} \right]^2 \intd A_s ^{R \mu}
	\right] 
	\nonumber
	\\ 
	&\le & 
	c_1 e^{-(R(t)+r_1-M)^2} t^{2(d-2)}
	E _{x} \left[
		\int _0 ^{t-s(t)-u} e^{A_s ^{(Q-1)\mu}} \intd A_s ^{R \mu}
	\right]
	\nonumber 
	\\ 
	&\le & 
	c_2 e^{c_3 t -(R(t)+r_1-M)^2} t^{2(d-2)}
	= 
	e^{-(R(t)+r_1-M)^2 + o(t^2)},  
	\label{eq:2ndorderpartless1} 
\eeq 

Next, suppose that $t-s(t)-u>1$. We set  
	\beqs
		&&  
		E _{x} \left[
			\int _0 ^{t-s(t)-u} e^{A_s ^{(Q-1)\mu}} 
			\bfE_{B_s} \left[ N_{t-s(t)-u-s} ^{C(R(t))} \right]^2 \intd A_s ^{R \mu}
		\right]	
		\\ 
		&= &
		E _{x} \left[ \int _0 ^{t-s(t)-u-1} \cdots \right]	
		+ 
		E _{x} \left[ \int ^{t-s(t)-u} _{t-s(t)-u-1} \cdots \right]	
		=
		I_1 + I_2 ,
	\eeqs 
	and $w_2=w_2 (s)=s(t)+u+s$, then $ 0 \le w_2 \le t-1$ for any $s \in [0,t-s(t)-u-1]$.
	Now (\ref{eq:20200106-1}) gives
	\beqs
	 &&
		\bfE_x \left[ N_{t-w_2} ^{C (R(t))} \right] 
		= 
		E_x \left[ e^{A_{t-w_2}^{(Q-1)\mu}} ; B_{t-w_2} \in C (R(t)) \right] 
		\\
		&=& 
		P_{x} \left( B_{t-w_2} \in C(R(t)) \right) + e^{-\lambda (t-w_2)} h(x) \int _{C(R(t))} h(y) \intd y 
		+ \int _{C(R(t))} q_{t-w_2} ^{(Q-1)\mu} (x,y) \intd y .
	\eeqs 
	Owing to Lemma \ref{lem:estP_0-1} (ii), (\ref{eq:c_*}) and Proposition \ref{prop:PropSN3-3} (i), 
	respectively, we have
	$$
		P_{x} \left( B_{t-w_2} \in C(R(t)) \right)
		\le 
		C e^{-\lambda (t-w_2)-\sqrt{-2\lambda}R(t)} t^{(d-2)/2} ,  
	$$ 
	$$
		e^{-\lambda (t-w_2)} h(x) \int _{C(R(t))} h(y) \intd y
		\le 
		C e^{-\lambda (t-w_2) - \sqrt{-2\lambda}R(t)} 
		R(t)^{(d-1)/2}
		\le 
		C' e^{-\lambda (t-w_2) - \sqrt{-2\lambda}R(t)} t^{(d-1)/2}
	$$
	and 
	$$ 
		\left|
			\int _{C(R(t))} q_{t-w_2} ^{(Q-1)\mu} (x,y) \intd y
		\right|	
		\le 
		C e^{-\lambda (t-w_2) - \sqrt{-2\lambda}R(t)} 
		\left( t^{(d-2)/2} \vee t^{(d-1)/2} \right) , 
	$$
	for any large $t$, $w_2 \in [0,t-1]$ and $x \in \ol{B_0(M)}$. 
	Hence
	\[
		\bfE _x \left[ N_{t-w_2} ^{C (R(t))} \right] ^2
		\le 
		C h^2(x) e^{-2\lambda (t-w_2) - 2 \sqrt{-2\lambda}R(t) +C_d(t)} ,
	\]
	which implies 
	\beqs
		I_1 
		&\le & 
		E _{x} \left[
			\int _0 ^{t-s(t)-u-1} e^{A_s ^{(Q-1)\mu}} 
			\bfE_{B_s} \left[ N_{t-w_2(s)} ^{C (R(t))} \right]^2 \intd A_s ^{R \mu}
		\right]	
		\\ 
		&\le & 
		C \| h \| _{\infty} ^2 e^{-2\lambda (t-s(t)-u) - 2\sqrt{-2\lambda} R(t) +C_d(t)}
		E _{x} \left[
			\int _0 ^{t-s(t)-u-1} e^{A_s ^{(Q-1)\mu} + 2 \lambda s} \intd A_s ^{R \mu}
		\right]	. 
	\eeqs 
	By the same argument as \cite[Proposition 3.3 (i)]{CS07}, we also have
	\[
		E _{x} \left[
			\int _0 ^{t-s(t)-u-1} e^{A_s ^{(Q-1)\mu} + 2 \lambda s} \intd A_s ^{R \mu}
		\right]
		\le 
		\sup _{x \in \rd} E _{x} \left[
			\int _0 ^{\infty } e^{A_s ^{(Q-1)\mu} + 2 \lambda s} \intd A_s ^{R \mu}
		\right]
		< \infty .
	\]
	Therefore, 
	\beq
		I_1 
		\le 
		C e^{-2\lambda (t-s(t)-u) - 2\sqrt{-2\lambda} R(t) + C_d (t)} .
	\label{eq:I1}
	\eeq 	
We set $w_3 = t-s(t)-u$, where $0 \le w_3-s \le 1$ for any $s \in [w_3-1, w_3] $. In the same way as 
(\ref{eq:20200226-4}), we have for any large $t$, 
\[
	\bfE_x \left[ N _{w_3-s} ^{C (R(t))} \right] 
	\le 
	C e^{-(R(t)+r_1-M)^2/2}t^{d-2} ,
\] 
and by the same as (\ref{eq:20201112-1}), 
$$
	E _{x} \left[ 
	\int _{w_3-1} ^{w_3} e^{A_s ^{(Q-1)\mu} } \intd A_s ^{R\mu}
	\right]
	\le 
	E _{x} \left[ 
		e^{A_{t} ^{(Q-1)\mu^+} } A^{R |\mu|} _{t}  
	\right] 	
	\le e^{c t} .
$$
Thus, 
\beqs
	I_2 
	&=& 
	E _{x} \left[ \int _{w_3-1} ^{w_3} e^{A_s ^{(Q-1)\mu} } 
		\bfE _{B_s} \left[ N_{w_3-s} ^{C(R(t))}\right]^2 \intd A_s ^{R \mu}
	\right] 
	\\ 
	&\le & 
	C e^{-(R(t) + r_1 -M)^2} t^{2(d-2)} 
		E _{x} \left[ \int _{w_3-1} ^{w_3} e^{A_s ^{(Q-1)\mu} } \intd A_s ^{R\mu}
		\right] 
	\\
	&\le & 
	C e^{-(R(t) + r_1 -M)^2 + c t} t^{2(d-2)} 
	= 
	e^{-(R(t) + r_1 -M)^2 + o(t^2)} .
\eeqs 
Combining this with (\ref{eq:I1}), we obtain 
\beq 
	E _{x} \left[ \int _0 ^{t-s(t)-u} e^{A_s ^{(Q-1)\mu} } 
		\bfE _{B_s} \left[ N_{t-s(t)-u-s} ^{C (R(t))} \right]^2 \intd A_s ^{R\mu}
	\right] 
	\le 
	C' e^{-2 \lambda (t-s(t)) -2\sqrt{-2\lambda} R(t) + C_d (t)} .
	\label{eq:2ndorderpart2less1}
\eeq
Consequently, (\ref{eq:2ndorderpartless1}) and (\ref{eq:2ndorderpart2less1}) yield that 
\beqs
	E _{x} \left[ \int _0 ^{t-s(t)-u} e^{A_s ^{(Q-1)\mu} } 
		\bfE _{B_s} \left[ N_{t-s(t)-u-s} ^{C (R(t))} \right]^2 \intd A_s ^{R\mu}
	\right] 
	&\le &
	C e^{-2\lambda (t-s(t)-u) - 2\sqrt{-2\lambda} R(t) + C_d (t)}
	\\ 
	&\le & 
	C' h(x) e^{-2\lambda (t-s(t)-u) - 2\sqrt{-2\lambda} R(t) + C_d (t)} , 
\eeqs 
for any large $t$, $u \in [0,t-s(t)]$ and $x \in \ol{B_0(M)}$. 
\epf
\begin{prop}\label{lem:2ndorder-estimate} 
We assume the same setting as in Proposition \ref{lem:orderof1moment}. 
\begin{enumerate}
\item[{\rm (i)}]
For any $\delta \in (0,\sqrt{-\lambda /2})$, there exist $T \ge 1$ and 
$\Theta (t)$ with $\Theta (t) \to 1$ as $t \to \infty$, such that for any $t \ge T$, 
	\beq
	\bfE _{x(t)} \left[ \left( N_{t-s(t)} ^{C(R(t))} \right)^2 \right] 
	\le 
	\Theta (t) h(x(t))e^{- 2 \lambda (t-s(t))- 2 \sqrt{-2\lambda}R(t) + C_d (t)} . 
	\label{eq:20210201-1}
	\eeq 
\item[{\rm (ii)}] Let $\delta \in [\sqrt{-\lambda/2} , \sqrt{-2\lambda})$. 
When $\delta = \sqrt{-\lambda /2}$, we set $a(t) = a_d(t)$ as in {\rm (\ref{eq:a(t)})} and assume the following: 
	\beq
		\lambda s(t)-\sqrt{-2\lambda} \gamma (t)
		\to 
		-\infty , \quad t \to \infty . 
		\label{eq:20210105-2}
	\eeq 
Then 
	\[
		\bfE _{x(t)} \left[ \left( N_{t-s(t)} ^{C(R(t))} \right)^2 \right] 
%		\sim  
%		\bfE _{x(t)} \left[ N_{t-s(t)} ^{C_r(R(t))} \right] 
		\sim 
		c_{*} h(x(t)) e^{-\lambda (t-s(t))-\sqrt{-2\lambda}R(t)} R(t)^{(d-1)/2}, 
		\quad t \to \infty . 
	\]
\end{enumerate}
\end{prop} 

\bpf 
	(i) If $\delta \in (0,\sqrt{-\lambda /2})$, 
	then $-\lambda (t-s(t)) - \sqrt{-2\lambda} R(t) > 0$, for large $t$. 
	Hence, Proposition \ref{lem:orderof1moment} and Lemma \ref{lem:2order} show that for any large $t$, 
	\beq
		\bfE _{x(t)} \left[ \left( N_{t-s(t)} ^{C(R(t))} \right)^2 \right] 
		&\le &
		c_* \theta _+ (t) h(x(t)) e^{-\lambda (t-s(t))-\sqrt{-2\lambda}R(t)}
		R(t)^{(d-1)/2}
		\nonumber
		\\ && \quad 
		+ 
		h(x(t)) e^{-2\lambda (t-s(t))-2\sqrt{-2\lambda}R(t) + C_d (t)} 
		\nonumber 
		\\ 
		&= & 
		\left( 
			\ds\frac{c_* \theta _+ (t) R(t)^{(d-1)/2}}
			{e^{-\lambda (t-s(t))-\sqrt{-2\lambda}R(t)+C_d(t)}}
			+ 
			1
		\right) 
		h(x(t)) e^{-2\lambda (t-s(t))-2\sqrt{-2\lambda}R(t) + C_d (t)}
		\nonumber
		\\
		&=:&
		\Theta (t) h(x(t)) e^{-2\lambda (t-s(t))-2\sqrt{-2\lambda}R(t) + C_d (t)} .
			\label{eq:20200924}
	\eeq 
	(ii) Let $\delta \in [\sqrt{-\lambda/2} , \sqrt{-2\lambda})$. 
		By the same way as {\rm (\ref{eq:20200924})}, 
		\beqs
			\bfE _{x(t)} \left[ \left( N_{t-s(t)} ^{C(R(t))} \right)^2 \right] 
			\hspace{-3mm} &\le& \hspace{-3mm}
			\left( 
				1 + 
				\ds\frac{e^{-\lambda (t-s(t))-\sqrt{-2\lambda}R(t)+C_d(t)}}{c_* \theta _+ (t) 
				R(t)^{(d-1)/2}}
			\right) 
			c_* \theta _+ (t)
			h(x(t)) e^{-\lambda (t-s(t))-\sqrt{-2\lambda}R(t)} 
			R(t)^{(d-1)/2}
			\\
			&=:& 
			c_* \Theta ' (t)
			h(x(t)) e^{-\lambda (t-s(t))-\sqrt{-2\lambda}R(t)} 
			R(t)^{(d-1)/2} . 
		\eeqs 
	Under the assumption, 
		\beq
			e^{-\lambda (t-s(t))-\sqrt{-2\lambda}R(t)+C_d(t)} R(t)^{-(d-1)/2} 
			\to 0 .
			\label{eq:criticalbehav}
		\eeq 
	Thus, we have $\Theta '(t) \to 1$. 
	By Proposition \ref{lem:orderof1moment} and Lemma \ref{lem:2order}, we also have the 
	corresponding lower estimate. 
	Thus,
	\[
		\bfE _{x(t)} \left[ \left( N_{t-s(t)} ^{C (R(t))} \right)^2 \right] 
		\sim  
		\bfE _{x(t)} \left[ N_{t-s(t)} ^{C(R(t))} \right] 
		\sim 
		c_{*} h(x(t)) e^{-\lambda (t-s(t))-\sqrt{-2\lambda}R(t)}
		R(t)^{(d-1)/2}, 
		\quad t \to \infty . 
	\]	
\epf 
\begin{rem}\label{rem:a(t)}
	In the case of $\delta = \sqrt{-\lambda/2}$, we need {\rm (\ref{eq:20210105-2})} 
	such that $\Theta '(t) \to 1$. 
	For example, if $s(t) =0$, then $\gamma (t)  = \log (t \vee 1)$;
	if $s(t)= \alpha \log (t \vee 1)$, $0 < \alpha < 1$, then $\gamma (t) \equiv \gamma \ge 0$; 
	if $s(t) \to \infty$ and $\gamma (t) \to -\infty$, then $\gamma (t) = o(s(t))$. 
\end{rem}

\section{Proofs}\label{sec:Bernoulli}
Let $\mu ^+$, $\mu^-$ be Kato class measures with compact supported in $\rd$. 
We assume Assumption \ref{ass}. 
Our proofs follow the same approach as \cite{Bocha2020} for the proofs of Theorems \ref{theorem1} 
and \ref{theorem2}. 
We also note that $R(t) =\delta t + a(t)$, $\delta \in (0 , \sqrt{-2\lambda})$ and $a(t)=o(t)$.  
\begin{prop}\label{prop:Bernoulli}
	Let $d \ge 1$ and $\delta \in [\sqrt{-\lambda/2}, \sqrt{-2\lambda})$. 
	When $\delta = \sqrt{-\lambda /2}$, we set $a(t) =a_d(t)$ as in {\rm (\ref{eq:a(t)})} and assume  
	that $\gamma (t)$ satisfies  {\rm (\ref{eq:20210105-2})}. 
	Then, there exist $\theta _i (t)$ {\rm(}$i=1,2,3,4${\rm)} with $\theta _i (t) \to 1$ as $t \to \infty$ 
	and $T > 0$ such that for all $t \ge T$, 
	\beq
		\bfP_{x(t)} \left( N_{t-s(t)} ^{C (R(t))} = 0 \right) 
		&\ge & 
		1- c_* h(x(t)) e^{-\lambda (t-s(t)) - \sqrt{-2\lambda} R(t)} 
		R (t)^{(d-1)/2} \theta _1 (t), 
		\label{eq:0-l}
		\\
		\bfP_{x(t)} \left( N_{t-s(t)} ^{C (R(t))} = 0 \right) 
		&\le & 
		1- c_* h(x(t)) e^{-\lambda (t-s(t)) - \sqrt{-2\lambda} R(t)} 
		R(t)^{(d-1)/2} \theta _2 (t), 
		\label{eq:0-u}
		\\
		\bfP_{x(t)} \left( N_{t-s(t)} ^{C (R(t))} = 1 \right) 
		&\ge & 
		c_* h(x(t)) e^{-\lambda (t-s(t)) - \sqrt{-2\lambda} R(t)} 
		R(t)^{(d-1)/2} \theta _3 (t), 
		\label{eq:1-l}
		\\ 
		\bfP_{x(t)} \left( N_{t-s(t)} ^{C (R(t))} = 1 \right) 
		&\le & 
		c_* h(x(t)) e^{-\lambda (t-s(t)) - \sqrt{-2\lambda} R(t)}
		R(t)^{(d-1)/2} \theta _4 (t), 
		\label{eq:1-u} 
		\\ 
		\bfP_{x(t)} \left( N_{t-s(t)} ^{C (R(t))} > 1 \right) 
		&\le & 
		h(x(t)) e^{-2 \lambda (t-s(t)) - 2\sqrt{-2\lambda} R(t) + C_d (t)} .
		\label{eq:otherwise}
	\eeq
	The constant $c_*$ is defined by {\rm(\ref{eq:c_*})} and 
	$C_d (t) = (d-1) \log (t \vee 1)$ for $d \ge 1$. 
\end{prop}
Before we prove Proposition \ref{prop:Bernoulli}, we note the following: 
Suppose that $Z$ is an $\n \cup \{0 \}$-valued random variable on some probability space 
$(\Omega , \calF , P)$. 
Then
	\beq
		P (Z>1) \le E \left[ Z^2  \right] - E [Z] \quad \text{and} \quad 
		\ds\frac{E \left[ Z\right]^2}{E \left[ Z^2 \right]} 
		\le 
		P \left( Z>0 \right) 
		\le 
		E \left[ Z \right] .
		\label{eq:Paley-Zygmund}
	\eeq 
\bpf 
	We first prove (\ref{eq:otherwise}). 
	Lemma \ref{lem:2order} and (\ref{eq:Paley-Zygmund}) give that for any large $t$, 
	\beqs
		\bfP _{x(t)} \left(N_{t-s(t)} ^{C (R(t))} > 1 \right) 
		&\le &
		\bfE _{x(t)} \left[ \left( N_{t-s(t)} ^{C (R(t))} \right) ^2 \right] - 
		\bfE _{x(t)} \left[ N_{t-s(t)} ^{C (R(t))} \right] 
		\\
		&\le & 
		h({x(t)}) e^{-2 \lambda (t-s(t)) - 2\sqrt{-2\lambda} R(t) + C_d (t))}. 
	\eeqs 
	We next prove (\ref{eq:0-l}). 
	On account of (\ref{eq:Paley-Zygmund}), we have
	\[ 
		\bfP _{x(t)} \left( N_{t-s(t)} ^{C (R(t))} = 0 \right) 
		=  
		1 - \bfP _{x(t)} \left( N_{t-s(t)} ^{C (R(t))} > 0 \right) 
		\ge 
		1 - \bfE _{x(t)} \left[ N_{t-s(t)} ^{C (R(t))} \right] . 
	\]
	Then it follows by Proposition \ref{lem:orderof1moment} that for any large $t$, 
	\[ 
		\bfP _{x(t)} \left( N_{t-s(t)} ^{C (R(t))} = 0 \right) 
		\ge 
		1 - \theta _+ (t) c_* h(x(t)) e^{-\lambda (t-s(t)) - \sqrt{-2\lambda}R(t)} 
		R(t)^{(d-1)/2}.
	\]
	
	We next prove (\ref{eq:0-u}). By (\ref{eq:Paley-Zygmund}), 
	\[ 
			\bfP _{x(t)} \left( N_{t-s(t)} ^{C (R(t))} = 0 \right) 
			=  
			1 - \bfP _{x(t)} \left( N_{t-s(t)} ^{C (R(t))} > 0 \right) 
			\le 
			1 -\ds\frac{ \bfE _{x(t)} \left[ N_{t-s(t)} ^{C (R(t))} \right]^2}
				{\bfE _{x(t)} \left[ \left( N_{t-s(t)} ^{C (R(t))} \right)^2 \right]} .
	\]
	According to Proposition \ref{lem:orderof1moment} and Lemma \ref{lem:2order}, we can choose 
	$\theta _- (t)$ and $\theta _+ ' (t)$ both converging to one as $t \to \infty$, 
	and for any large $t$, 
	\beqs
		&& 
		\ds\frac{ \bfE _{x(t)} \left[ N_{t-s(t)} ^{C(R(t))} \right]^2}
				{\bfE _{x(t)} \left[ \left( N_{t-s(t)} ^{C (R(t))} \right)^2 \right]}
		\ge 
		\ds\frac{ \bfE _{x(t)} \left[ N_{t-s(t)} ^{C (R(t))} \right]^2 }
				{\bfE _{x(t)} \left[ N_{t-s(t)} ^{C (R(t)) } \right] 
			+C h (x(t)) e^{-2\lambda (t-s(t)) - 2 \sqrt{-2\lambda} R(t) + C_d (t)} }
		\\
		&\ge &
		\ds\frac{ \left( c_* \theta _- (t) h(x(t)) 
		e^{-\lambda (t-s(t)) - \sqrt{-2\lambda}R(t)} R(t)^{(d-1)/2} \right)^2}
		{c_* \theta _+' (t) h(x(t)) e^{-\lambda (t-s(t)) - \sqrt{-2\lambda}R(t)}
		R(t)^{(d-1)/2}
		+ Ch (x(t)) e^{-2\lambda (t-s(t)) - 2 \sqrt{-2\lambda} R(t) + C_d(t)} }
		\\
		&=:& 
		c_* \theta _2 (t) h(x(t)) e^{-\lambda (t-s(t)) - \sqrt{-2\lambda}R(t)}
		R(t)^{(d-1)/2}, 
	\eeqs
	where by (\ref{eq:criticalbehav}), 
	\beq 
		\theta _2 (t)
		= 
		\ds\frac{ \theta _-^2 (t) }
		{\theta _+ ' (t) 
		+ 
		c_*^{-1} C
		e^{-\lambda (t-s(t)) - \sqrt{-2\lambda}R(t) + C_d(t)} R(t)^{-(d-1)/2}
		} \to 1 , \quad t \to \infty .
		\label{eq:bottleneck}
	\eeq 
	Therefore, 
	\[	
		\bfP _{x(t)} \left( N_{t-s(t)} ^{C (R(t))} = 0 \right) 
		\le 
		1 - 
		c_* \theta _2 (t) h(x(t)) e^{-\lambda (t-s) - \sqrt{-2\lambda}R(t)} 
		R(t)^{(d-1)/2}.
	\] 
	Since $\theta _2(t) \to 1$ as $t \to \infty$, we have (\ref{eq:0-u}). 
	
	We finally prove (\ref{eq:1-l}) and (\ref{eq:1-u}). 
	By (\ref{eq:0-l}), (\ref{eq:0-u}) and (\ref{eq:otherwise}), 
	\beqs
		&& 
		\bfP _{x(t)} \left( N_{t-s(t)} ^{C (R(t))} = 1 \right) 
		= 
		1- \bfP _{x(t)} \left( N_{t-s(t)} ^{C (R(t))} > 1 \right) - \bfP _{x(t)} \left( N_{t-s(t)} ^{C (R(t))} = 0 \right) 
		\\ 
		&\ge &
		1 - h(x(t)) e^{-2\lambda (t-s(t)) -2\sqrt{-2\lambda} R(t) + C_d(t)}- 
		\left( 1 - c_* h(x(t)) e^{-\lambda (t-s(t)) - \sqrt{-2\lambda} R(t) } \theta _2 (t)
		R(t)^{(d-1)/2} \right)
		\\
		&=&
		\left( 
			\theta _2 (t) 
			- 
			\ds\frac{e^{-\lambda (t-s(t)) -\sqrt{-2\lambda} R(t) + C_d (t)} }
			{c_* R(t)^{(d-1)/2}}
		\right)c_* h(x(t)) e^{-\lambda (t-s(t)) - \sqrt{-2\lambda} R(t) } 
		R(t)^{(d-1)/2}
		\\
		&=:& 
		\theta _3 (t) c_* h(x(t)) e^{-\lambda (t-s) - \sqrt{-2\lambda} R(t) }
		R(t)^{(d-1)/2}. 
	\eeqs 
	Similarly to (\ref{eq:bottleneck}), $\theta _3 (t) \to 1$ as $t \to \infty$. 
	By (\ref{eq:Paley-Zygmund}) and Proposition \ref{lem:orderof1moment}, 
	there exists $\theta _4 (t)$ such that $\theta _4 (t) \to 1$ as $t \to \infty$ and
	\beqs 
		\bfP _{x(t)} \left( N_{t-s(t)} ^{C (R(t))} = 1 \right) 
		&\le & 
		\bfP _{x(t)} \left( N_{t-s(t)} ^{C (R(t))} > 0 \right) 
		\le 
		\bfE _{x(t)} \left[N_{t-s(t)} ^{C (R(t))} \right] 
		\\
		&\le & 
		c_* h(x(t)) e^{-\lambda (t-s(t)) - \sqrt{-2\lambda} R(t) } R(t)^{(d-1)/2} \theta _4 (t) .
	\eeqs 
\epf 
The rest of this section will be devoted to the proofs of our main theorems. 
\subsection{Proof of Theorem \ref{theorem1}}
\bpf
	In Proposition \ref{prop:Bernoulli}, we choose $x(t) \equiv x$ and $s(t) \equiv 0$. 
	In particular, for $\delta = \sqrt{-\lambda /2}$, 
	we set $a(t)=a_d(t)$ as (\ref{eq:a(t)}) and assume $0 \le \gamma (t) \to \infty$.  
	Then $a_d (t)$ satisfies (\ref{eq:20210105-2}) in Proposition \ref{prop:Bernoulli}. 
	By (\ref{eq:0-l}) and (\ref{eq:0-u}), 
	there exist $T > 0$ and $\theta _i(t)$, $i=1,2$ such that for any $t \ge T$, 
	\[  
		c_* h(x) e^{-\lambda t - \sqrt{-2 \lambda} R(t)} R(t)^{(d-1)/2} \theta _2 (t)
		\le 
		\bfP _x \left( N_t ^{C (R(t))} > 0 \right) 
		\le 
		c_* h(x) e^{-\lambda t - \sqrt{-2 \lambda} R(t)} R(t)^{(d-1)/2} \theta _1 (t) . 
	\]  
	Hence, we have the conclusion. 
\epf 

\subsection{Proof of Theorem \ref{theorem2}}
Let us denote by $Z_t (u)$ the descendants of $u$ at $t$ and 
$Z_t ^A (u)$ the subset of $Z_t (u)$ whose particles are contained in domain $A$. 
In addition, we use $N_t (u)$ and $N_t ^A (u)$ as the size of $Z_t (u)$ and $Z_t ^A (u)$, respectively. 
We see from the Markov property and variance formula that
\beq
	&&
	\bfE _{x} \left[ \left. 
		\left( N_t ^{C(R(t))} - \bfE_{\bfB_{s(t)}} \left[ N_{t-s(t)} ^{C(R(t))} \right] \right)^2
		\, \right| \calG _{s(t)} 
	\right] 
	\nonumber
	\\
	&=& 
	\sum _{u \in Z_{s(t)}} 
	\left\{
		\bfE _{\bfB_{s(t)}^u} \left[ \left( N_{t-s(t)}^{C(R(t))} (u) \right)^2 \right]
		-
		\bfE _{\bfB_{s(t)}^u} \left[ N_{t-s(t)}^{C(R(t))} (u) \right]^2
	\right\} 
	\nonumber
	\\
	&\le &
	\sum _{u \in Z_{s(t)}} 
		\bfE _{\bfB_{s(t)}^u} \left[ \left( N_{t-s(t)}^{C(R(t))} (u) \right)^2 \right]	.
	\label{eq:varianceformula} 
\eeq 
\begin{proof}[Proof of Theorem \ref{theorem2}]
Since $M_t \to M_{\infty}$, $\bfP _x $-a.s. and
\beqs
		&& 
	\left|
		e^{\lambda t + \sqrt{-2\lambda} R(t)} R(t)^{-(d-1)/2}
		N_t ^{C (R(t))}- c_* M_{\infty}
	\right|
	\\ 
	& \le & 
	e^{\lambda t + \sqrt{-2\lambda} R(t)} R(t)^{-(d-1)/2}
	\left|
		N_t ^{C (R(t))} - \bfE_{\bfB_{s(t)}} \left[ N_{t-s(t)} ^{C(R(t))} \right]
	\right|
	\\ 
	&& \ 
	+
	\left|
		e ^{\lambda t + \sqrt{-2\lambda} R(t)} R(t)^{-(d-1)/2}
		\bfE _{\bfB _{s(t)}} \left[ N_{t-s(t)} ^{C(R(t))} \right] 
		- c_* M_{s(t)}
	\right|	
	+ 
	c_* \left| M_{s(t)} - M_{\infty} \right|, 
\eeqs 
we need to show the following:  
\beq
	&& 
	e^{\lambda t + \sqrt{-2\lambda} R(t)} R(t)^{-(d-1)/2}
	\left| N_t ^{C (R(t))} - \bfE_{\bfB_{s(t)}} \left[ N_{t-s(t)} ^{C(R(t))} \right] \right| 
	\to 0, \quad \text{ in  probability } \bfP _x ,
	\label{eq:convinprob2}
	\\ 
	&& 
	\left| 
		e ^{\lambda t + \sqrt{-2\lambda} R(t)} R(t)^{-(d-1)/2}
		\bfE _{\bfB _{s(t)}} \left[ N_{t-s(t)} ^{C(R(t))} \right] 
		- 
		c_* M_{s(t)}
	\right| \to 0, \quad \bfP _x\text{-a.s.}.
	\label{eq:convinprob1}
\eeq 
We first prove (\ref{eq:convinprob2}). 
Fix $\varepsilon _0 > 0$. For any $\varepsilon > 0$, 
\beq
\begin{split}
	&
	\bfP _x \left( e^{\lambda t + \sqrt{-2\lambda} R(t)} 
		R(t)^{-(d-1)/2}
		\left| N_t ^{C(R(t))} - \bfE_{\bfB_{s(t)}} \left[ N_{t-s(t)} ^{C(R(t))} \right] \right| > \varepsilon 
	\right)
	\\ 
	= \, \, & 
	\bfP _x \left( e^{\lambda t + \sqrt{-2\lambda} R(t)} R(t)^{-(d-1)/2}
			\left| N_t ^{C(R(t))} - \bfE_{\bfB_{s(t)}} \left[ N_{t-s(t)} ^{C(R(t))} \right] \right| > \varepsilon , \ 
			L_{s(t)} \le \left( \sqrt{\frac{-\lambda}{2}} + \varepsilon _0 \right) s(t)
	\right)
	\\ 
	&
	+ 
	\bfP _x \left( e^{\lambda t + \sqrt{-2\lambda} R(t)} R(t)^{-(d-1)/2}
			\left| N_t ^{C(R(t))} - \bfE_{\bfB_{s(t)}} \left[ N_{t-s(t)} ^{C(R(t))} \right] \right| > \varepsilon, \ 
			L_{s(t)} > \left( \sqrt{\frac{-\lambda}{2}} + \varepsilon _0 \right) s(t)
	\right) .
\end{split} 
\label{eq:20200205-5} 
\eeq
Since $L_t / t \to \sqrt{-\lambda /2}$, $\bfP _x$-a.s. by (\ref{Ltlim}), 
the second term of (\ref{eq:20200205-5}) converges to zero. 
We set $b(t) = (\sqrt{-\lambda /2} + \varepsilon_0)s(t)$, where $s(t) \to \infty$. 
Let $\delta \in (0,\sqrt{-\lambda/2})$. 
Then the Chebyshev inequality and (\ref{eq:varianceformula}) yield 
\beq
	&& 
	\bfP _x \left( 
		e^{\lambda t + \sqrt{-2\lambda} R(t)} R(t)^{-(d-1)/2}
		\left| N_t ^{C(R(t))} - \bfE_{\bfB_{s(t)}} \left[ N_{t-s(t)} ^{C(R(t))} \right]  \right| > \varepsilon, \, 
		L_{s(t)} \le b(t)
	\right)
	\nonumber 
	\\ 
	&=& 
	\bfE _x \left[ 
	\bfP _{x} \left( \left. 
		e^{\lambda t + \sqrt{-2\lambda} R(t)} R(t)^{-(d-1)/2}
		\left| N_t ^{C(R(t))} - \bfE_{\bfB_{s(t)}} \left[ N_{t-s(t)} ^{C(R(t))} \right]  \right| > \varepsilon
		\, \right| \calG _{s(t)}
	\right)
	; L_{s(t)} \le b(t)
	\right] 
	\nonumber 
	\\
	&\le & 
	\ds\frac{e^{2\lambda t + 2\sqrt{-2\lambda} R(t)}R(t)^{1-d}}{\varepsilon ^2}
	\bfE _x \left[ 
	\bfE _{x} \left[ \left. 
		\left( N_t ^{C(R(t))} - \bfE_{\bfB_{s(t)}} \left[ N_{t-s(t)} ^{C(R(t))} \right] \right)^2
		\, \right| \calG _{s(t)}
	\right]
	; L_{s(t)} \le b(t)
	\right] 
	\nonumber
	\\
	&\le&
	\ds\frac{e^{2\lambda t + 2\sqrt{-2\lambda} R(t)} R(t)^{1-d}}{\varepsilon ^2}
	\bfE _x \left[ 
	\sum _{u \in Z_{s(t)}} 
		\bfE _{\bfB_{s(t)}^u} \left[ \left( N_{t-s(t)}^{C(R(t))} (u) \right)^2 \right]
	; L_{s(t)} \le b(t)
	\right] . 
	\label{eq:20200130-0}
\eeq
We note that $|\bfB _{s(t)} ^u| \le b (t)$ for any $u \in Z_{s(t)}$ on the event $\{ L_{s(t)} \le b(t) \}$. 
According to Proposition \ref{lem:2ndorder-estimate} (i), we can take non-random $T>0$ so large that 
for all $t >T$, the second moment in (\ref{eq:20200130-0}) is bounded by (\ref{eq:20210201-1}) 
uniformly on the event $\{ L_{s(t)} \le b(t) \}$. That is, 
\beqs
	(\ref{eq:20200130-0})
	&\le & 
	\ds\frac{\Theta (t) e^{\lambda s(t) + C_d(t)} R(t)^{1-d}}{\varepsilon ^2}
	\bfE _x \left[ 
			e^{\lambda s(t)}
			\sum _{u \in Z_{s(t)}} h \left( \bfB _{s(t)} ^u \right) 			
		; L_{s(t)} \le b(t)
	\right]	
	\\
	&\le & 
	\ds\frac{\delta^{1-d} \Theta (t) e^{\lambda s(t)}}{\varepsilon ^2}
	\bfE _x \left[ M_{s(t)}\right]	
	= 
	\ds\frac{\delta ^{1-d}\Theta (t) e^{\lambda s(t)}}{\varepsilon ^2} h(x) 
	\to 0 , \quad t \to \infty. 
\eeqs
We now consider the critical case for (\ref{eq:convinprob2}). 
Let us suppose that $\gamma (t) \to - \infty$ and $\gamma(t) = o(\log t)$. 
We take $s(t) = O (\log t)$, which satisfies both (\ref{eq:20210105-2}) and the condition as mentioned at 
the beginning of Section \ref{sec3}. 
By the same argument as (\ref{eq:20200130-0}) and Proposition \ref{lem:2ndorder-estimate} (ii), 
\beqs
	&& 
	\bfP _x \left( 
		e^{\lambda t + \sqrt{-2\lambda} R(t)} R(t)^{-(d-1)/2}
		\left| N_t ^{C(R(t))} - \bfE_{\bfB_{s(t)}} \left[ N_{t-s(t)} ^{C(R(t))} \right]  \right| > \varepsilon, \, 
		L_{s(t)} \le b(t)
	\right)
	\\
	&\le &
	\ds\frac{e^{2\lambda t + 2\sqrt{-2\lambda} R(t)} R(t)^{1-d}}{\varepsilon ^2}
	\bfE _x \left[ 
	\sum _{u \in Z_{s(t)}} 
		\bfE _{\bfB_{s(t)}^u} \left[ \left( N_{t-s(t)}^{C(R(t))} (u) \right)^2 \right]
	; L_{s(t)} \le b(t)
	\right] 
	\\ 
	&\le &
	\ds\frac{e^{2\lambda t + 2\sqrt{-2\lambda} R(t)} R(t)^{1-d}}{\varepsilon ^2}
	C e^{-\lambda t -\sqrt{-2\lambda} R(t)} R(t)^{(d-1)/2}
		\bfE _x \left[ 
				e^{\lambda s(t)}
				\sum _{u \in Z_{s(t)}} h \left( \bfB _{s(t)} ^u \right) 			
			; L_{s(t)} \le b(t)
		\right]	
	\\
	&\le & 
	\dfrac{C}{\varepsilon^2} 
	e^{\lambda t + \sqrt{-2\lambda} R(t)} R(t)^{(1-d)/2} \bfE _x \left[ M_{s(t)} \right]
	= 
	\dfrac{C}{\varepsilon^2} h(x)
	e^{\lambda t + \sqrt{-2\lambda} R(t)} R(t)^{(1-d)/2} . 
\eeqs
Since 
\[
	e^{\lambda t + \sqrt{-2\lambda} R(t)} R(t)^{(1-d)/2}
	\sim 
	\left( \dfrac{-\lambda}{2} \right)^{(1-d)/4} e^{\sqrt{-2 \lambda} a(t) + \{(1-d)/2\} \log t } 
	=
	\left( \dfrac{-\lambda}{2} \right)^{(1-d)/4} e^{\sqrt{-2 \lambda}\gamma (t)},
\]
the second term of (\ref{eq:20200205-5}) also converges to zero. 

We next prove (\ref{eq:convinprob1}). 
By abuse of notation, we use 
\[ 
	A (t)= e ^{\lambda t + \sqrt{-2\lambda} R(t)} R(t)^{-(d-1)/2}
	\bfE _{\bfB_{s(t)}} \left[ N_{t-s(t)} ^{C(R(t))}\right] .
\]
Then we see that for any $\delta \in (0, \sqrt{-2\lambda})$, $A(t) \rightarrow c_* M_\infty$, $\bfP _x$-a.s. 
In fact, by Proposition \ref{lem:orderof1moment}, there exists non-random $T' >T$ such that for any $t>T'$, 
we have uniformly on $\{ L_{s(t)} \le b(t) \}$, 
\beqs
	A(t) &=& 
	e ^{\lambda t + \sqrt{-2\lambda}R(t)} R(t)^{-(d-1)/2}
	\sum _{u \in Z_{s(t)}}
	\bfE _{\bfB_{s(t)}^u} \left[ N_{t-s(t)} ^{C(R(t))} (u)\right]
	\\
	&\le & 
	e ^{\lambda t + \sqrt{-2\lambda}R(t)}
	R(t)^{-(d-1)/2}
	\sum _{u \in Z_{s(t)}}
	\theta_+ (t) c_* h(\bfB^u_{s(t)})
	e ^{-\lambda (t-s(t)) - \sqrt{-2\lambda} R(t)} 
	R(t)^{(d-1)/2}
	\\
	&=& 
	c_* \theta_+ (t) M_{s(t)} .
\eeqs  
We also have $A (t) \ge c_* \theta_- (t) M_{s(t)}$ so that 
$A(t) \rightarrow c_* M_\infty$, $\bfP _x$-a.s. 
We thus have (\ref{eq:convinprob1}). 
\epf

{\subsection{Proof of Theorem \ref{theorem3}}
Let $d=1,2$ and $R(t) = \sqrt{-\lambda /2} t + a(t)$, where $a(t)=a_d (t)$ given by (\ref{eq:a(t)}). 
We  can choose $s(t)$ which satisfies both (\ref{eq:20210105-2}) and 
the condition as mentioned at the beginning of Section. 
If $\gamma (t) \ge 0$, then we can take such an $s(t)$ independently of $\gamma (t)$, 
for example $s(t) = \alpha \log (t \vee 1)$ for any $\alpha \in (0,1)$.  
Hence, we can use Proposition \ref{prop:Bernoulli} to prove Theorem \ref{theorem3}. 

For a fixed $\varepsilon > 0$, we use $b(t)= (\sqrt{-\lambda /2} + \varepsilon ) s(t)$
as in (\ref{eq:20200205-5}).
\begin{lem}\label{lem:20210308}
	There exists a non-random $T$ and $\kappa _1 (t)$, $\kappa_2 (t)$ such that for all $t \ge T$,  
	\beq 
		\begin{split}
		\exp \left( -\kappa _1 (t) c_* (i) M_{s(t)} e^{-\sqrt{-2\lambda} a(t)} 
		R(t)^{(d-1)/2} \right)
		&\le 
		\prod _{u \in Z_{s(t)}} \bfP _{\bfB _{s(t)} ^u} 
		\left(
			N_{t-s(t)} ^{C^i (R(t))} = 0
		\right) 	
		\\ 
		& \le 
		\exp \left( -\kappa _2 (t) c_* (i) M_{s(t)} e^{-\sqrt{-2\lambda} a(t)} 
		R(t)^{(d-1)/2}
		\right) 
		\end{split}
		\label{eq:basiceq} 
	\eeq
	uniformly on the event $\{ L_{s(t)} \le b(t) \}$. 
	Here $\kappa _1 (t), \kappa _2 (t) \to 1$ as $ t\to \infty$. 
\end{lem}
		By (\ref{Ltlim}), $\indi _{[0,b(t)]} (L_{s(t)}) = \indi _{[0, \sqrt{-\lambda/2} + \varepsilon]} (L_{s(t)}/s(t))
		\to 1$, $\bfP _x$-a.s. 
		Since 
		$	e ^{- \sqrt{-2\lambda} a(t)} R(t) ^{(d-1)/2}
			\sim 
			e^{- \sqrt{-2\lambda} \gamma (t)}
		$, we see from (\ref{eq:basiceq}) that 
		\beq
		\begin{split}
			&
			\lim _{t \to \infty} \indi _{[0,b(t)]} \left( L_{s(t)} \right) \prod _{i=1} ^n \prod _{v \in Z_{s(t)}} 
			\bfP _{\bfB _{s(t)}^v} \left( N_{t-s(t)} ^{C^i (R(t))} (v) = 0 \right)
			\\ 
			&=
			\begin{cases}
			\exp \left( - \wt{c} M_{\infty} e^{-\sqrt{-2\lambda} \gamma} \right)  , 
			& \text{if } \gamma (t) \to \gamma < \infty , \\ 
			1, & \text{if } \gamma (t) \to \infty ,
			\end{cases}		
			\quad \bfP _x\text{-a.s.} 
		\end{split} 
			\label{eq:P1-1} 
		\eeq 
\bpf[Proof of Lemma \ref{lem:20210308}] 
	By Proposition \ref{prop:Bernoulli}, there exists a non-random $T>0$ such that for any $t>T$ and $x(t)$ 
	with $|x(t)| \le b(t)$, 
	\[
		\bfP _{x(t)} \left( N_{t-s(t)} ^{C^ i(R(t))} = 0 \right)
		\le 
		1- 
		c_*(i) h \left( x(t) \right) e^{\lambda s(t) - \sqrt{-2\lambda} a(t)} 
		R(t)^{(d-1)/2}
		\theta _2 (t) .
	\]
	Since $|\bfB _{s(t)} ^u| \le b(t)$ for all $u \in Z_{s(t)}$ on the event $\{ L_{s(t)} \le b(t) \}$ 
	and $1-x \le e^{-x}$ for all $x \in \R$, 
	\beqs 
		\prod _{u \in Z_{s(t)}} 
			\bfP _{\bfB _{s(t)}^u} \left( N_{t-s(t)} ^{C^ i(R(t))} = 0 \right)
		&\le& 
		\exp \left\{
		- \sum _{u \in Z_{s(t)}} c_* (i) h \left( \bfB _{s(t)} ^u \right) 
		e^{\lambda s(t) - \sqrt{-2\lambda} a(t)} 
		R(t)^{(d-1)/2}
		\theta _2 (t) 
		\right\} 
		\\ 
		&=& 
		\exp \left(  - c_* (i) \theta _2 (t) M_{s(t)} e^{- \sqrt{-2\lambda} a(t)} 
		R(t)^{(d-1)/2}
		\right) .
	\eeqs 
	For a fixed $x^* \in (0,1)$, $\log (1-x) \ge \frac{\log (1-x^*)}{x^*} x$ for any $x  \in (0,x^*) $. 
	Then by Proposition \ref{prop:Bernoulli}, for any $t>T$, 
	\beqs
		&& 
		\prod _{u \in Z_{s(t)}} 
			\bfP _{\bfB _{s(t)}^u} \left( N_{t-s(t)} ^{C^i (R(t))} (u) = 0 \right)
		\ge 
		\prod _{u \in Z_{s(t)}} 
		\left(
		1- 
		c_* (i) h \left( \bfB _{s(t)} ^u \right) e^{\lambda s(t) - \sqrt{-2\lambda} a(t)} 
		R(t)^{(d-1)/2}
		\theta _1 (t) 
		\right)
		\\
		&=& 
		\exp \left\{ 
			\sum _{u \in Z_{s(t)}} 
			\log \left( 1- c_* (i) h \left( \bfB _{s(t)} ^u \right) e^{\lambda s(t) - \sqrt{-2\lambda} a(t)} 
			R(t)^{(d-1)/2}
			\theta _1 (t) \right) 
		\right\} 
		\\
		&\ge & 
		\exp \left\{ 
			\sum _{u \in Z_{s(t)}} 
			\ds\frac{\log \left( 1- c_* (i) \| h \|_{\infty} e^{\lambda s(t) - \sqrt{-2\lambda} a(t)} 
				R(t)^{(d-1)/2} \theta _1 (t) \right)}
				{c_* (i) \| h \|_{\infty} e^{\lambda s(t) - \sqrt{-2\lambda} a(t)} 
				R(t)^{(d-1)/2} \theta _1 (t)} 
		\right. 
		\\
		&& \qquad \times 
		\left.
			c_* (i) h \left( \bfB _{s(t)} ^u \right) e^{\lambda s(t) - \sqrt{-2\lambda} a(t)} 
			R(t)^{(d-1)/2}\theta _1 (t)
		\right\} 
		\\
		&=& 
		\exp \left\{
			\ds\frac{\log \left( 1- c_* (i) \| h \|_{\infty} e^{\lambda s(t) - \sqrt{-2\lambda}a(t)} 
			R(t)^{(d-1)/2} \theta _1 (t) \right)}
				{c_* (i) \| h \|_{\infty} e^{\lambda s(t) - \sqrt{-2\lambda} a(t)} R(t)^{(d-1)/2}
				\theta _1 (t)} 
				c_* (i) \theta _1 (t) M_{s(t)}
				e^{- \sqrt{-2\lambda} a(t)} 
				R(t)^{(d-1)/2}
		\right\} .
	\eeqs  
	We note that 
	\[
		e^{\lambda s(t) - \sqrt{-2\lambda} a(t)} R(t)^{(d-1)/2}
		\sim  
		\left( \dfrac{-\lambda}{2} \right) ^{(d-1)/4} 
		e^{\lambda s(t) - \sqrt{-2\lambda} \gamma (t)} 
		\to 0, 
	\]
	by (\ref{eq:20210105-2}). 
	Since $\frac{\log (1-x)}{x} \to -1$ as $x \downarrow 0$, we have
	\[
		\kappa _1 (t) := -
			\ds\frac{\log \left( 1- c_* (i) \| h \|_{\infty} e^{\lambda s(t) - \sqrt{-2\lambda} a(t)} 
			R(t)^{(d-1)/2} \theta _1 (t) \right)}
				{c_* (i) \| h \|_{\infty} e^{\lambda s(t)- \sqrt{-2\lambda} a(t)} 
				R(t)^{(d-1)/2} \theta _1 (t)} \theta _1 (t)
		\to 1, \quad t \to \infty ,
	\]
	so that (\ref{eq:basiceq}) follows. 
\epf

\bpf[Proof of Theorem \ref{theorem3}]
	For any fixed $\varepsilon >0$, we set $b(t) = (\sqrt{-\lambda /2} + \varepsilon ) s(t)$ and  
	$a (t) = a_d (t)$ as (\ref{eq:a(t)}). 
	Then
	\beq
	\begin{split}
		\bfP _x \left( \bigcap _{i=1}^n \left\{ N_t ^{C^i (R(t))} = k _i \right\} \right)
		&= 
		\bfP _x \left( \bigcap _{i=1}^n \left\{ N_t ^{C^i (R(t))} = k _i \right\} , L_{s(t)} \le b(t) \right)
		\\ 
		& \qquad 
		+
		\bfP _x \left( \bigcap _{i=1}^n \left\{ N_t ^{C^i (R(t))} = k _i \right\} , L_{s(t)} > b(t) \right) 
	\end{split}
	\label{eq:20200908-1}
	\eeq
	and the second term converges to zero by (\ref{Ltlim}). 
	We here consider the limit of the first term in (\ref{eq:20200908-1}). 
	If $k=0$, then 
	\beqs
		&& 
		\bfP _x \left( \bigcap _{i=1}^n \left\{ N_t ^{C^i (R(t))} = 0 \right\} , L_{s(t)} \le b(t) \right)
		\\
		&=&
		\bfE _x \left[ \bfP _{\bfB _{s(t)}} 
			\left( \bigcap _{i=1}^n \left\{ N_{t-s(t)} ^{C^i (R(t))} = 0 \right\} \right) ; L_{s(t)} \le b(t) \right]	
		\\
		&=& 
		\bfE _x \left[ \prod _{i=1} ^n \prod _{u \in Z_{s(t)}} \bfP _{\bfB _{s(t)}^u} 
			\left( N_{t-s(t)} ^{C^i (R(t))} = 0 \right) ; L_{s(t)} \le b(t) \right]	.
	\eeqs 
	We see from (\ref{eq:basiceq}), (\ref{eq:P1-1}) and the bounded convergence theorem that, 
	\beqs 
		&& 
		\liminf _{t \to \infty} 
		\bfP _x \left( \bigcap _{i=1}^n \left\{ N_t ^{C^i (R(t))} = 0 \right\} , L_{s(t)} \le b(t) \right)
		\\ 
		&\ge&
		\liminf _{t \to \infty} 
		\bfE _x \left[
			\prod _{i=1}^n
			\exp \left( - \kappa _1 (t) c_* (i) M_{s(t)} e^{- \sqrt{-2\lambda} a(t)} 
			R(t)^{(d-1)/2} \right) ; L_{s(t)} \le b(t)
		\right]
		\\
		&\ge &
		\lim _{t \to \infty} 
		\bfE _x \left[
			\prod _{i=1}^n
			\exp \left( - \kappa _1 (t) c_* (i) M_{s(t)} e^{- \sqrt{-2\lambda} a(t)} 
			R(t)^{(d-1)/2} \right) 
		\right] 
		\\
		&=& 
		\begin{cases}
		\bfE _x \left[ \exp \left( -\wt{c} M_{\infty} e^{- \sqrt{-2\lambda} \gamma} \right) \right] , 
		& \text{if } \gamma(t) \to \gamma < \infty , \\ 
		1 , & \text{if } \gamma (t) \to \infty .
		\end{cases} 
	\eeqs 	
	Here, the second inequality above is adapted from $\bfP _x (L_{s(t)} > b(t)) \to 0$ as $t \to \infty$. 
	Hence the limit is one when $\gamma(t) \to \infty$. 
	In the case of $\gamma(t) \to \gamma$,
	\beqs 
		&& 
		\limsup _{t \to \infty} 
		\bfP _x \left( \bigcap _{i=1}^n \left\{ N_t ^{C^i (R(t))} = 0 \right\}, L_{s(t)} \le b(t) \right)
		\\ 
		& \le & 
		\limsup _{t \to \infty} 
		\bfE _x \left[
			\prod _{i=1}^n
			e^{- \kappa _2 (t) c_* (i) M_{s(t)} \sqrt{-2 \lambda }a(t)} R(t)^{(d-1)/2} ; L_{s(t)} \le b(t)
		\right]
		\\
		&\le &  
		\lim _{t \to \infty} 
		\bfE _x \left[
			\prod _{i=1}^n
			e^{- \kappa _2 (t) c_* (i) M_{s(t)} - \sqrt{-2 \lambda }a(t)} R(t) ^{(d-1)/2} 
		\right]
		= 
		\bfE _x \left[ \exp \left( - \wt{c} M_{\infty} e^{- \sqrt{-2\lambda} \gamma} \right) \right]  .
	\eeqs 
	Therefore, our claim holds for $k=0$. 
	
	We suppose that $k \ge 1$. 
	Let $[m]$ denote a subset of $Z_{s(t)}$ with size $m \ge 1$ and $[0]=\emptyset$. 
	Under the law $\bfP _{\bfB_{s(t)}}$, 
	the event $\bigcap _{i=1} ^n \{ N_{t-s(t)} ^{C^i (R(t))} = k _i \} $ is divided into 
	$B$ and its complementary event, where $B$ is defined by the following: 
	$B(0)=\Omega$, $\sum _{j=1} ^0 = \emptyset$ and inductively, 
	\beqs
		&& 
		B (i) := 
		\Biggl\{
			\omega \in B(i-1) \ \left| \
			{}^{\exists} [k_i] \subset Z_{s(t)} (\omega) \setminus \sum _{j=1}^{i-1}[k_j] \right.
						\text{ such that } 
			\\ 
			&&
			N_{t-s(t)}^{C^i (R(t))} (u) (\omega) = 1, \ \text{for all } u \in [k_i] 
						\text{ and }
						N_{t-s(t)}^{C^i (R(t))} (v)(\omega) = 0, \ \text{for all } v \in Z_{s(t)} (\omega)
						\setminus [k_i] 
		\Biggr\} , 
	\eeqs 
	then $B := B(n)$. 
	Here, $N_{t-s(t)} ^{C^i (R(t)) } (u)$ is the number of particles which are descendants of $u$ and 
	lie on $C^i (R(t))$ at $t-s(t)$. 
	By the strong Markov property,  
	\beqs
		&& 
		\bfP _x \left( \bigcap _{i=1}^n \left\{ N_t ^{C^i (R(t))} = k _i \right\} , L_{s(t)} \le b(t) \right)
		=
		\bfE _x \left[ \bfP _{\bfB _{s(t)}} 
			\left( \bigcap _{i=1}^n \left\{ N_{t-s(t)} ^{C^i (R(t))} = k _i \right\} \right) ; L_{s(t)} \le b(t) \right]
		\\
		&=&
		\bfE _x \left[ \bfP _{\bfB _{s(t)}} \left( B \right) ; L_{s(t)} \le b(t) \right]
		+ 
		\bfE _x \left[ \bfP _{\bfB _{s(t)}} 
			\left( B^c \cap \bigcap _{i=1}^n \left\{ N_{t-s(t)} ^{C^i (R(t))} = k _i \right\} \right) 
			; L_{s(t)} \le b(t) \right]
		=:
		E_1 + E_2 .
	\eeqs
	Since 
	\[ 
		B^c \cap \bigcap _{i=1}^n \left\{ N_{t-s(t)} ^{C^i (R(t))} = k _i \right\} 
		\subset 
		\left\{
			\omega \ \left| \ {}^\exists u \in Z_{s(t)} (\omega) \text{ such that } 
			N_{t-s(t)} ^{C^i(R(t))} (u) (\omega) >1, \right. 
			\text{for some $i$}
		\right\}
	\]
	Proposition \ref{prop:Bernoulli} implies that  
	\beqs
		E_2 
		&\le & 
		\sum _{i=1} ^n
		\bfE _x \left[ \sum _{u \in Z_{s(t)}} \bfP _{\bfB _{s(t)}^u} \left(  N_{t-s(t)} ^{C^i(R(t))} > 1 \right) ; 
		L_{s(t)} \le b(t) \right] 	
		\\
		&\le &
		n
		\bfE _x \left[ \sum _{u \in Z_{s(t)}} h \left( \bfB_{s(t)} ^u \right) e^{-2\lambda (t-s(t)) 
			- 2 \sqrt{-2\lambda}R(t) + C_d (t)} ; 
		L_{s(t)} \le b(t) \right] 	
		\\
		&=&
		n
		e^{\lambda s(t) - 2\sqrt{-2\lambda} a(t) +C_d(t)}
		\bfE _x \left[ e^{\lambda s(t)} \sum _{u \in Z_{s(t)}} h \left( \bfB_{s(t)} ^u \right) ; 
		L_{s(t)} \le b(t) \right] 
		\\
		&= &
		ne^{\lambda s(t)- 2 \sqrt{-2\lambda} a(t)+C_d(t)} \bfE _x \left[ M_{s(t)} \right] 
		=   
		nh(x) e ^{\lambda s(t)- 2 \sqrt{-2\lambda} \gamma (t)} \to 0, \quad t \to \infty .
	\eeqs 
	In what follows, we compute the limit of $E_1$. 
	Let $(k)$ denote a $k$-sequence $(u_1 , \dots , u_k)$, $u_i \in Z_{s(t)}$ and $u_i \neq u_j$.  
	We distinguish between $k$-permutations. 
	We admit the following two claims: 
			\beq
			E_1 =
				\ds\frac{1}{k_1 ! k_2 ! \dots k_n !}
				\bfE _x \left[  
				\prod_{j=1} ^n  
				\prod _{v \in Z_{s(t)}} \bfP _{\bfB _{s(t)}^v} \left( N_{t-s(t)} ^{C^j(R(t))} (v) = 0 \right) 	
				\right. 
				\nonumber
				\\  
				\left. \quad \times 
				\sum _{(k) \subset Z_{s(t)}} 
				\prod_{i=1} ^n \prod _{u \in k_i(k)} 
				\ds\frac{\bfP _{\bfB _{s(t)}^{u}} \left( N_{t-s(t)} ^{C^i(R(t))} (u) =1 \right)}
				{\bfP _{\bfB _{s(t)}^{u}} \left( N_{t-s(t)} ^{C^i(R(t))} (u) =0 \right)}
				; L_{s(t)} \le b(t) 
				\right] 
			\label{eq:20200310-3}
			\eeq 
		and 
			\beq
			\begin{split}
				& \lim _{t \to \infty} \indi _{[0,b(t)]} \left( L_{s(t)} \right)
				\sum _{(k) \subset Z_{s(t)}} \prod _{i=1} ^n \prod _{u \in k_i(k)}
					\ds\frac{\bfP _{\bfB _{s(t)}^{u}} \left( N_{t-s(t)} ^{C_r(R(t))} (u) =1 \right) }
						{\bfP _{\bfB _{s(t)}^{u}} \left( N_{t-s(t)} ^{C_r(R(t))} (u) = 0 \right)} 
				\\ 
				&=
				\begin{cases} 
				\left(  M_{\infty} e^{-\sqrt{-2\lambda} \gamma} \right)^k 
				\ds\prod _{i=1} ^n c_* (i) ^{k_i} , & \text{if } \gamma (t) \to \gamma < \infty, \\
				0, & \text{if } \gamma (t) \to \infty,
				\end{cases} 
				\quad \bfP_x \text{-a.s.},
			\end{split}
			\label{eq:P1-2}
			\eeq		
		which will be proved in Lemmas \ref{lem:4-15}, \ref{lem:4-16} below.  
		As is clear from the definition of $E_1$, we can use the bounded convergence theorem. 
		By (\ref{eq:P1-1}), (\ref{eq:20200310-3}), (\ref{eq:P1-2}), 
		\[
			\lim _{t \to \infty} E_1 
			= 
			\begin{cases} 
			\bfE _x \left[ \left(  M_{\infty} e^{-\sqrt{-2\lambda} \gamma }  \right)^k 
			\exp \left( -\wt{c} M_{\infty} e^{-\sqrt{-2\lambda} \gamma} \right) \right] 
			\ds\prod _{i=1} ^n \frac{c_* (i) ^{k_i}}{k_i !}, & \text{if } \gamma (t) \to \gamma < \infty , \\ 
			0, & \text{if } \gamma (t) \to \infty , 
			\end{cases} 
		\]
		for any $k \ge 1$, which completes the proof. 
	\epf
	\begin{lem}\label{lem:4-15}
		For $k \ge 1$, {\rm (\ref{eq:20200310-3})} holds. 
	\end{lem}
	\bpf
		Each $k$-sequence $(k)$ has $k_i$-subsequences ($i=1,\dots, n$) 
		such that $(k) = ( (k_1), (k_2) , \dots , (k_n))$. 
		We will use $k_i (k)$ to represent the $k_i$-subsequences. 
		In addition, we abuse notation slightly  
		$Z_{s(t)} \setminus k_i (k) = \{ u \in Z_{s(t)} \mid u \not\in k_i (k) \}$. 
		We define the equivalence relation $\sim _k$ on $\{ (k) \mid (k) \subset Z_{s(t)} \}$ 
		if and only if $k_i (k)_1 = k_i (k)_2$ as two sets, for all $i=1,\dots ,n$, for $(k)_1$, $(k)_2 \subset Z_{s(t)}$.  
		Then we denote by $Z_{s(t)} \slash \sim _k $ the equivalence class $\{ (k) \mid (k) \subset Z_{s(t)} \}$ 
		by $\sim _k$. 
		Hence, 
		\beqs 
			B
			&=&  
			\bigcup _{(k) \in Z_{s(t)} \slash \sim _k } \bigcap _{i=1} ^n 
				\left\{ 
					N_{t-s(t)} ^{C^i(R(t))} (u) =1, {}^{\forall} u \in k_i (k) \text{ and } 
					N_{t-s(t)} ^{C^i(R(t))} (v) = 0 , {}^{\forall} v \in Z_{s(t)} \setminus k_i (k)
				\right\} 
			\\ 
			&=:& 
			\sum _{(k)\in Z_{s(t)} \slash \sim _k } \bigcap _{i=1} ^n B ^i _{(k)} ,
		\eeqs 
		We note that $\bigcap _{i=1} ^n B_{(k)} ^i$ and $\bigcap _{i=1} ^n B_{(k)'} ^i$ are disjoint events 
		if $(k)$ and $(k)'$ are different sets. 
		Then, for a sequence $(k) \subset Z_{s(t)}$, 
		$\bigcap _{i=1} ^n B_{(k)} ^i$ has $k_1 ! \cdots k_n !$ multiples. 
		Let $\sum _{(k) \subset Z_{s(t)}}$ denote the summation over all $k$-sequences of $Z_{s(t)}$. 
		Thus, 
		\beqs
			E_1
			&=& 
			\bfE _x \Biggl[
			\bfP _{\bfB _{s(t)}} 
			\left( \bigcup _{(k) \in Z_{s(t)} \slash \sim _k } \bigcap _{i=1} ^n B^i _{(k)} \right) 
				; L_{s(t)} \le b(t) 
			\Biggr]
			\\
			&=& 
			\bfE _x \left[  \sum _{(k) \in Z_{s(t)} \slash \sim _k} 
			\bfP _{\bfB _{s(t)}} 
				\left( \bigcap _{i=1} ^n B_{(k)} ^i \right) ; L_{s(t)} \le b(t) 
			\right] 
			\\
			&=& 
			\frac{1}{k_1 ! k_2 ! \dots k_n !}
			\bfE _x \left[  
				\sum _{(k) \subset Z_{s(t)}} 
				\bfP _{\bfB _{s(t)}} 
				\left( \bigcap _{i=1} ^n B_{(k)} ^i \right) ; L_{s(t)} \le b(t) 
			\right]
		\eeqs  
		and 
		\beq 
			&& 
			\sum _{(k) \subset Z_{s(t)}} 
						\bfP _{\bfB _{s(t)}} 
						\left( \bigcap _{i=1} ^n B_{(k)} ^i \right)
			\nonumber
			\\ 
			&=& 
			\sum _{(k) \subset Z_{s(t)}} \prod_{i=1} ^n \left( 
			\prod _{u \in k_i(k)} \bfP _{\bfB _{s(t)}^{u}} \left( N_{t-s(t)} ^{C^i(R(t))} (u) =1 \right) 
			\prod _{v \in Z_{s(t)} \setminus k_i (k)} 
			\bfP _{\bfB _{s(t)}^v} \left( N_{t-s(t)} ^{C^i(R(t))} (v) = 0 \right) \right) 
			\nonumber 
			\\
			&=& 
			\sum _{(k) \subset Z_{s(t)}} \prod_{i=1} ^n 
			\left( 
				\prod _{v \in Z_{s(t)}} 
				\bfP _{\bfB _{s(t)}^v} \left( N_{t-s(t)} ^{C^i(R(t))} (v) = 0 \right) 
				\prod _{u \in k_i(k)} 
				\ds\frac{\bfP _{\bfB _{s(t)}^{u}} \left( N_{t-s(t)} ^{C^i(R(t))} (u) =1 \right)}
				{\bfP _{\bfB _{s(t)}^{u}} \left( N_{t-s(t)} ^{C^i(R(t))} (u) =0 \right)} 		
			\right) 
			\nonumber
			\\
			&=& 
			\prod_{j=1} ^n  
			\prod _{v \in Z_{s(t)}} \bfP _{\bfB _{s(t)}^v} \left( N_{t-s(t)} ^{C^j(R(t))} (v) = 0 \right) 	
			\sum _{(k) \subset Z_{s(t)}} 
			\prod_{i=1} ^n \prod _{u \in k_i(k)} 
			\ds\frac{\bfP _{\bfB _{s(t)}^{u}} \left( N_{t-s(t)} ^{C^i(R(t))} (u) =1 \right)}
			{\bfP _{\bfB _{s(t)}^{u}} \left( N_{t-s(t)} ^{C^i(R(t))} (u) =0 \right)} . 
			\nonumber
		\eeq
		Hence we have (\ref{eq:20200310-3}).
	\epf 
	
	\begin{lem}\label{lem:4-16}
		For $k \ge 1$, {\rm (\ref{eq:P1-2})} holds. 
	\end{lem}
	\bpf 
	Let $\gamma (t) \to \gamma \in [0,\infty)$. 
	We first show 
		\beq
		\begin{split}
		& 
			\liminf _{t \to \infty} \indi _{[0,b(t)]} \left( L_{s(t)} \right)
			\sum _{(k) \subset Z_{s(t)}} \prod _{i=1} ^n \prod _{u \in k_i (k)}  
				\ds\frac{\bfP _{\bfB _{s(t)}^{u}} \left( N_{t-s(t)} ^{C^i(R(t))} (u) =1 \right) }
					{\bfP _{\bfB _{s(t)}^{u}} \left( N_{t-s(t)} ^{C^i(R(t))} (u) = 0 \right)} 
		\\
		&
			\ge 
				\left( e^{-\sqrt{-2\lambda} \gamma} M_{\infty}\right)^k
			\prod _{i=1} ^n c_* (i) ^{k_i} , \quad \bfP_x \text{-a.s.} 
		\end{split}
			\label{eq:20200701-1}
		\eeq	
		We note that 
		\beq
			e^{-\lambda t - \sqrt{-2\lambda} R(t)} R(t) ^{(d-1)/2} 
			= 
			e^{-\sqrt{-2\lambda} \gamma (t) + o(1)} \quad \text{ as } t \to \infty .
			\label{eq:20210401-1}
		\eeq
		By Proposition \ref{prop:Bernoulli},  
		\beq
			&& 
			\sum _{(k) \subset Z_{s(t)}} \prod _{i=1}^n \prod _{u \in k_i(k)}  
				\ds\frac{\bfP _{\bfB _{s(t)}^{u}} \left( N_{t-s(t)} ^{C^i(R(t))} (u) =1 \right) }
					{\bfP _{\bfB _{s(t)}^{u}} \left( N_{t-s(t)} ^{C^i(R(t))} (u) = 0 \right)} 
			\nonumber
			\\
			&\ge &
			\sum _{(k) \subset Z_{s(t)}} \prod _{i=1} ^n \prod _{u \in k_i(k)}  
				\bfP _{\bfB _{s(t)}^{u}} \left( N_{t-s(t)} ^{C^i (R(t))} (u) =1 \right) 
			\nonumber 
			\\ 
			&\ge & 
			\sum _{(k) \subset Z_{s(t)}} \prod _{i=1} ^n \prod _{u \in k_i(k)}  
				c_* (i)h \left( \bfB _{s(t)}^{u} \right)	
				e^{-\lambda (t-s(t)) - \sqrt{-2\lambda} R(t)} 
				R(t)^{(d-1)/2}
				\theta _3 (t)
			\nonumber
			\\ 
			&= &
			\left( e^{-\sqrt{-2\lambda} \gamma (t) + o(1)} \theta _3 (t) \right)^k 
			\sum _{(k) \subset Z_{s(t)}} \prod _{i=1} ^n \prod _{u \in k_i(k)}   
				c_* (i)h \left( \bfB _{s(t)}^{u} \right) e^{\lambda s(t)} .
			\label{eq:20200311-2}
		\eeq
		For $k \ge 2$, we set $(Z_{s(t)})^k = \{ (u_1 , \dots , u_k) \mid u_i \in Z_{s(t)} , 1 \le i \le k \}$ 
		and 
		\[
			\Lambda = \left\{ (u_1 , \dots , u_k) \mid {}^{\exists} \{ p,q \} \text{ such that } u_p = u _q  \right\} 
			\subset \left( Z_{s(t)} \right)^k .
		\]
		For $U \in (Z_{s(t)})^k$, we write $U=(u_1,\dots,u_k)$. 
		Since the totality of the $k$-permutations of $Z_{s(t)}$ is in the one-to-one correspondence to
		$(Z_{s(t)})^k \setminus \Lambda$, 
		\beq
			&& 
			\sum _{(k) \subset Z_{s(t)}} \prod _{i=1} ^n \prod _{u \in k_i(k)}  
			c_* (i) h \left( \bfB _{s(t)}^{u} \right)	e^{\lambda s(t)}
			\nonumber 
			\\
			&=&
			\sum _{U \in (Z_{s(t)})^k} \prod _{i=1} ^n \prod _{u \in k_i(U)}   
				c_* (i)h \left( \bfB _{s(t)}^{u} \right) e^{\lambda s(t)} 
			-
			\sum _{V \in \Lambda} \prod _{i=1} ^n 
					\prod _{u \in k_i(V)}  
						c_* (i)h \left( \bfB _{s(t)}^{u} \right) e^{\lambda s(t)} , 
			\label{eq:20200311-4} 
		\eeq 
		where 
		$
			k_i (U) = ( u_{\ol{k}_{i-1} + 1} , \dots, u_{\ol{k}_i} )  
		$ 
		is a subsequence of $U$ and $\ol{k}_0=0$, $\ol{k} _i = k_1 + \dots + k_{i}$ for $1 \le i \le n$. 
		Inductively, 
		\beq
			&&
			\sum _{U \in (Z_{s(t)})^k} \prod _{i=1} ^n 
			\prod _{u \in k_i (U)}  
				c_* (i)h \left( \bfB _{s(t)}^{u} \right) e^{\lambda s(t)} 
			\nonumber
			\\ 
			&=& 
			\sum _{v \in Z_{s(t)}} \sum _{U \in (Z_{s(t)})^{k-1}} 
			\left( 
				\prod _{i=1} ^{n-1} \prod _{j=\ol{k}_{i-1}+1} ^{\ol{k}_i}
				c_* (i)h \left( \bfB _{s(t)}^{u_j} \right) e^{\lambda s(t)}
			\right) 
			\nonumber
			\\ 
			&& \times 
			\left( 
				c_*(n) h \left( \bfB _{s(t)} ^{v} \right) e^{\lambda s(t)}
				\prod _{j=\ol{k}_{n-1}+1} ^{k -1} c_* (n) h \left( \bfB _{s(t)}^{u_j} \right) e^{\lambda s(t)} 
			\right) 
			\nonumber
			\\
			&=& 
			c_* (n)M_{s(t)} 
			\sum _{U \in (Z_{s(t)})^{k-1}} 
			\left( 
				\prod _{i=1} ^{n-1} \prod _{j=\ol{k}_{i-1}+1} ^{\ol{k}_i}
				c_* (i)h \left( \bfB _{s(t)}^{u_j} \right) e^{\lambda s(t)}
			\right) 
			\left( 
				\prod _{j=\ol{k}_{n-1}+1} ^{k -1} c_* (n) h \left( \bfB _{s(t)}^{u_j} \right) e^{\lambda s(t)} 
			\right) 
			\nonumber 
			\\
			&=& 
			M_{s(t)} ^k \prod _{i=1} ^n c_* ^{k_i} (i) .
			\label{eq:20200311-3}
		\eeq
		We divide $\Lambda $ into $\Lambda_{p,q} (\cdot )$ as follows: 
		\beqs 
			\Lambda 
			&=& 
			\sum _{1 \le p < q \le k} \left\{ (u_1 , \dots, u_k) \in \left( Z_{s(t)} \right)^k \mid u_p = u_q \right\}
			\\ 
			&=& 
			\sum_{u \in Z_{s(t)}}
			\sum _{1 \le p < q \le k} \left\{ (u_1 , \dots, u_k) \in \left( Z_{s(t)} \right)^k \mid u_p = u_q = u \right\}
			\\ 
			&=:& 
			\sum_{u \in Z_{s(t)}}\sum _{1 \le p < q \le k} \Lambda _{p,q} (u) . 
		\eeqs 
		For $V \in \Lambda _{p,q} (\cdot )$, we define 
		\[
			V | _{p,q} = \left\{ \left( u_1 , \dots , u_{p-1}, u_{p+1}, \dots , u_{p-1}, u_{q+1}, \dots u_k \right)
			\mid u_i \in Z_{s(t)}  \right\}
			= \left( Z_{s(t)} \right) ^{k-2} .
		\] 
		Let $C= \ds\max_{1 \le i \le n} c_* (i)$. 
		Then the second term of (\ref{eq:20200311-4}) is bounded below by 
		\beq
			&&
			-C \sum _{V \in \Lambda} \prod _{i=1} ^n 
				\prod _{v \in k_i (V)} h \left( \bfB _{s(t)}^{v} \right) e^{\lambda s(t)}
			=
			-C	\sum_{u \in Z_{s(t)}}\sum _{1 \le p < q \le k} \sum _{V \in \Lambda _{p,q} (u)}
				\prod _{i=1} ^n 
				\prod _{v \in k_i (V)} h \left( \bfB _{s(t)}^{v} \right) e^{\lambda s(t)}
			\nonumber
			\\
			&=&
			- C \sum_{u \in Z_{s(t)}} \left( h \left( \bfB _{s(t)}^{u} \right) e^{\lambda s(t)} \right) ^2
			\sum _{1 \le p < q \le k} \sum _{V \in \Lambda _{p,q} (u) }
			\prod _{v \in V |_{p,q}} h \left( \bfB _{s(t)}^{v} \right)	e^{\lambda s(t)}
			\nonumber 
			\\
			&=& 
			- C \sum_{u \in Z_{s(t)}} \left( h \left( \bfB _{s(t)}^{u} \right) e^{\lambda s(t)} \right) ^2
			\sum _{1 \le p < q \le k} \sum _{V \in (Z_{s(t)})^{k-2} }
			\prod _{v \in V} h \left( \bfB _{s(t)}^{v} \right)	e^{\lambda s(t)} .
			\label{eq:20210106-1}
		\eeq
		Similarly to (\ref{eq:20200311-3}), 
		\beq
			(\ref{eq:20210106-1}) 
			=
			- C \binom{k}{2} M_{s(t)} ^{k-2}
			\sum_{u \in Z_{s(t)}} \left( h \left( \bfB _{s(t)}^{u} \right) e^{\lambda s(t)} \right) ^2
			\ge 
			-C \binom{k}{2}
			M_{s(t)} ^{k-1} \| h \| _{\infty} e^{\lambda s(t)}.
			\label{eq:20200311-5} 
		\eeq
		Since we see from (\ref{eq:20200311-4}), (\ref{eq:20200311-3}) and (\ref{eq:20200311-5}) that 
		$$
		\sum _{(k) \subset Z_{s(t)}} \prod _{i=1} ^n \prod _{u \in k_i(k)}  
			c_* (i) ^{k_i} h \left( \bfB _{s(t)}^{u} \right)	e^{\lambda s(t)}
		\ge  
			M_{s(t)} ^k \prod _{i=1} ^n c_* (i) ^{k_i} 
			- C' M_{s(t)} ^{k-1} e^{\lambda s(t)}, 
		$$
		it follows from (\ref{eq:20200311-2}) that, 
		\beqs
			&& 
			\liminf _{t \to \infty} \indi _{[0,b(t)]} \left( L_{s(t)} \right)
			\sum _{(k) \subset Z_{s(t)}} \prod _{i=1} ^n \prod _{u \in k_i (k)}  
				\ds\frac{\bfP _{\bfB _{s(t)}^{u}} \left( N_{t-s(t)} ^{C^i(R(t))} (u) =1 \right) }
					{\bfP _{\bfB _{s(t)}^{u}} \left( N_{t-s(t)} ^{C^i(R(t))} (u) = 0 \right)} 
			\nonumber 
			\\ 
			&\ge & 
			\liminf _{t \to \infty} \indi _{[0,b(t)]} \left( L_{s(t)} \right) 
			\left( e^{-\sqrt{-2\lambda} \gamma (t) + o(1)} \theta _3 (t) \right)^k  
			\left( M_{s(t)} ^k \prod _{i=1} ^n c_*(i) ^{k_i} -
			C' M_{s(t)} ^{k-1} e^{\lambda s(t)} \right) 
			\nonumber
			\\
			&=&
				\left(  M_{\infty} e^{-\sqrt{-2\lambda} \gamma} \right)^k
			\prod _{i=1} ^n c_* (i) ^{k_i} , \quad \bfP_x \text{-a.s.} 
		\eeqs
		
		Assume that $k=1$, that is, $k_i=1$ for some $i$ and the others are zero. 
		By (\ref{eq:20200311-2}), 
		\beqs 
		&&
			\sum _{(k) \subset Z_{s(t)}} \prod _{i=1}^n \prod _{u \in k_i(k)}  
				\ds\frac{\bfP _{\bfB _{s(t)}^{u}} \left( N_{t-s(t)} ^{C^i(R(t))} (u) =1 \right) }
				{\bfP _{\bfB _{s(t)}^{u}} \left( N_{t-s(t)} ^{C^i(R(t))} (u) = 0 \right)} 
		\ge 
				e^{-\sqrt{-2\lambda} \gamma (t) + o(1)} 
			\theta _3 (t) c_* (i) \sum _{u \in Z_{s(t)}} h \left( \bfB _{s(t)}^{u} \right)	e^{\lambda s(t)}
			\\ 
		&=& 
			e^{-\sqrt{-2\lambda} \gamma (t) + o(1)} 
			\theta _3 (t) c_* (i) M_{s(t)}
		\to 
				c_* (i)  M_{\infty} e^{-\sqrt{-2\lambda} \gamma} .
		\eeqs 
		Hence (\ref{eq:20200701-1}) also holds for $k=1$. 
	
		We next show
		\beq
		\begin{split}
			& \limsup _{t \to \infty} \indi _{[0,b(t)]} 
			\left( L_{s(t)} \right) \sum _{(k) \subset Z_{s(t)}} \prod _{i=1} ^n \prod _{u \in k_i (k)}  
				\ds\frac{\bfP _{\bfB _{s(t)}^{u}} \left( N_{t-s(t)} ^{C^i(R(t))} (u) =1 \right) }
					{\bfP _{\bfB _{s(t)}^{u}} \left( N_{t-s(t)} ^{C^i(R(t))} (u) = 0 \right)}
			\\ 
			&
			\le
			\left(  M_{\infty} e^{-\sqrt{-2\lambda} \gamma} \right) ^k 
			\prod_{i=1} ^n c_* (i) ^{k_i} , \quad \bfP_x \text{-a.s.}
		\end{split}
			\label{eq:20200908-2}
		\eeq 
		By Proposition \ref{prop:Bernoulli} and (\ref{eq:20210401-1}),   
		\beq
			&& 
			\sum _{(k) \subset Z_{s(t)}} \prod _{i=1} ^n \prod _{u \in k_i(k)} 
				\ds\frac{\bfP _{\bfB _{s(t)}^{u}} \left( N_{t-s(t)} ^{C^i(R(t))} (u) =1 \right) }
					{\bfP _{\bfB _{s(t)}^{u}} \left( N_{t-s(t)} ^{C^i(R(t))} (u) = 0 \right)}
			\nonumber
			\\
			&\le &
			\sum _{(k) \subset Z_{s(t)}} \prod _{i=1} ^n \prod _{u \in k_i(k)}
				\ds\frac{c_* (i) h \left( \bfB_{s(t)}^{u} \right) 
					e^{-\lambda (t-s(t)) - \sqrt{-2\lambda}R(t)}
					R(t)^{(d-1)/2}
					\theta _4 (t)}
					{1-c_* (i) h \left( \bfB_{s(t)}^{u} \right) 
					e^{-\lambda (t-s(t)) - \sqrt{-2\lambda}R(t)} R(t)^{(d-1)/2} \theta _1 (t)}
			\nonumber 
			\\
			&= &
			\sum _{(k) \subset Z_{s(t)}} \prod _{i=1} ^n \prod _{u \in k_i(k)} 
					\ds\frac{c_* (i) e^{-\sqrt{-2 \lambda} \gamma (t) + o(1) } \theta _4 (t)}
						{1-c_* (i) \| h\|_\infty 
						e^{\lambda s(t) -\sqrt{-2 \lambda} \gamma (t) + o(1) } \theta _1 (t)}
				h \left( \bfB_{s(t)}^{u} \right) e^{\lambda s(t)} .
			\label{eq:20200311-1}
		\eeq
		By (\ref{eq:20200311-4}) and (\ref{eq:20200311-3}),  for any $k \ge 2$, 
		\beqs
		(\ref{eq:20200311-1})
		&\le& 
				\sum _{U \in (Z_{s(t)})^k} \prod _{i=1} ^n \prod _{u \in k_i(U)} 
					\ds\frac{c_* (i) e^{-\sqrt{-2 \lambda} \gamma (t) + o(1)} \theta _4 (t)}
						{1-c_* (i) \| h\|_\infty 
						e^{\lambda s(t) -\sqrt{-2 \lambda} \gamma (t) + o(1)} \theta _1 (t)}
				h \left( \bfB_{s(t)}^{u} \right) e^{\lambda s(t)} 
		\\ 
		&=&
		M_{s(t)} ^k
			\prod _{i=1} ^n 
			\left( 
				\ds\frac{c_* (i) e^{-\sqrt{-2 \lambda} \gamma (t) + o(1)} \theta _4 (t)}
				{1-c_* (i) \| h\|_\infty 
				e^{\lambda s(t) -\sqrt{-2 \lambda} \gamma (t) + o(1)} \theta _1 (t)}
			\right)^{k_i}.
		\eeqs 
		Combining this with (\ref{eq:20200311-1}), we have
		\beqs
			&& 
			\limsup _{t \to \infty} \indi _{[0,b(t)]} 
			\left( L_{s(t)} \right) 
			\sum _{(k) \subset Z_{s(t)}} \prod _{i=1} ^n \prod _{u \in k_i(k)} 
				\ds\frac{\bfP _{\bfB _{s(t)}^{u}} \left( N_{t-s(t)} ^{C^i(R(t))} (u) =1 \right) }
					{\bfP _{\bfB _{s(t)}^{u}} \left( N_{t-s(t)} ^{C^i(R(t))} (u) = 0 \right)}
			\nonumber
			\\
			&\le & 
			\limsup _{t \to \infty} \indi _{[0,b(t)]} \left( L_{s(t)} \right) M_{s(t)} ^k
			\prod _{i=1} ^n 
			\left( 
				\ds\frac{c_* (i) e^{-\sqrt{-2 \lambda} \gamma (t) + o(1)} \theta _4 (t)}
				{1-c_* (i) \| h\|_\infty 
				e^{\lambda s(t) -\sqrt{-2 \lambda} \gamma (t) + o(1)} \theta _1 (t)}
			\right)^{k_i}
			\nonumber
			\\ 
			&=& 
			\left(  M_{\infty} e ^{-\sqrt{-2\lambda}} \right)^k 
			\prod_{i=1} ^n c_* (i) ^{k_i} , \quad \bfP_x \text{-a.s.}
		\eeqs
		Assume that $k=1$, that is, $k_i = 1$ some $i$ and the others are zero. 
		By (\ref{eq:20200311-1}), 
		\beqs 
			\sum _{(k) \subset Z_{s(t)}} \prod _{i=1} ^n \prod _{u \in k_i(k)}
				\ds\frac{\bfP _{\bfB _{s(t)}^{u}} \left( N_{t-s(t)} ^{C^i(R(t))} (u) =1 \right) }
					{\bfP _{\bfB _{s(t)}^{u}} \left( N_{t-s(t)} ^{C^i(R(t))} (u) = 0 \right)}
			&\le& 
				\ds\frac{c_* (i) e^{-\sqrt{-2 \lambda} \gamma (t) + o(1)} \theta _4 (t)}
				{1-c_* (i) \| h\|_\infty 
				e^{\lambda s(t) -\sqrt{-2 \lambda} \gamma (t) + o(1)} \theta _1 (t)}
			M_{s(t)} 
			\\ 
			&\to & 
			c_* (i) M_{\infty} e^{-\sqrt{-2\lambda} \gamma} .
		\eeqs 
		Hence, (\ref{eq:20200908-2}) also holds for $k=1$. 
		From (\ref{eq:20200701-1}) and (\ref{eq:20200908-2}), we obtain the first part of (\ref{eq:P1-2}). 
		
		When $\gamma (t) \to \infty$, by the same way as (\ref{eq:20200311-1}), 
		\beqs
		&&
			\limsup _{t \to \infty} 
			\sum _{(k) \subset Z_{s(t)}} \prod _{i=1} ^n \prod _{u \in k_i(k)} 
				\ds\frac{\bfP _{\bfB _{s(t)}^{u}} \left( N_{t-s(t)} ^{C^i(R(t))} (u) =1 \right) }
					{\bfP _{\bfB _{s(t)}^{u}} \left( N_{t-s(t)} ^{C^i(R(t))} (u) = 0 \right)}		
		\\ 
		&\le& 
			\limsup _{t \to \infty} 
			M_{s(t)} ^k
			\prod _{i=1} ^n 
			\left( 
				\ds\frac{c_* (i) e^{-\sqrt{-2 \lambda} \gamma (t) + o(1)} \theta _4 (t)}
				{1-c_* (i) \| h\|_\infty 
				e^{\lambda s(t) -\sqrt{-2 \lambda} \gamma (t) + o(1)} \theta _1 (t)}
			\right)^{k_i}
			= 0.
		\eeqs 
		Thus we have (\ref{eq:P1-2}). 
	\epf

\renewcommand{\abstractname}{Acknowledgements}
\begin{abstract}
	The author is deeply grateful to Professor Yuichi SHIOZAWA for a careful reading of the manuscript and 
	many helpful comments. 
	The author also would like to thank Christopher B. Prowant for carefully proofreading the manuscript. 
\end{abstract}

\end{document}